\documentclass[11pt]{article}

\usepackage{amscd}
\usepackage{latexsym}
\usepackage{amsmath,amsfonts,amssymb}
\usepackage{t-angles}
  \providecommand{\Pairing}[1]{\Put(10,20)[cc]{#1}\ev}

\numberwithin{equation}{section}
\newcommand{\sstyle}{\scriptstyle}

\newcommand{\Label}[1]{\label{#1}\ \\ \makebox[0mm][r]
  {\fontsize{9}{11pt}\fontseries{bx}\fontshape{n}
  \selectfont #1\hspace*{8mm}}\!\!}
 \newcommand{\mathlabel}[1]{\label{#1}\makebox[0mm][r]
  {\fontsize{9}{11pt}\fontseries{bx}\fontshape{n}
  \selectfont #1\hspace*{8mm}}\!\!}
\renewcommand{\Label}{\label}        
\renewcommand{\mathlabel}{\label}    

\newcommand{\crossbox}[4]{\raisebox{3pt}{\def\arraystretch{.5}%
\begin{tabular}{|c|c|}\hline\raisebox{-.7pt}{$\sstyle\!\!\!#1\!\!\!$}&%
\raisebox{-.7pt}{$\sstyle\!\!\!#2\!\!\!$}\\%
\hline\raisebox{-.7pt}{$\sstyle\!\!\!#3\!\!\!$}&%
\raisebox{-.7pt}{$\sstyle\!\!\!#4\!\!\!$}\\\hline\end{tabular}}}
\newcommand{\itembox}[4]{\item[\crossbox{#1}{#2}{#3}{#4}\ :]}

\newtheorem{theorem}{{\bf{T}}{\fontsize{7pt}{11}\fontshape{n}
\fontseries{bx}\selectfont{$\!\!\!$HEOREM}}}[section]
\newtheorem{proposition}[theorem]{{\bf{P}}{\fontsize{7pt}{11}\fontshape{n}
\fontseries{bx}\selectfont{$\!\!\!$ROPOSITION}}}
\newtheorem{corollary}[theorem]{{\bf{C}}{\fontsize{7pt}{11}\fontshape{n}
\fontseries{bx}\selectfont{$\!\!\!$OROLLARY}}}
\newtheorem{lemma}[theorem]{{\bf{L}}{\fontsize{7pt}{11}\fontshape{n}
\fontseries{bx}\selectfont{$\!\!\!$EMMA}}}
\newtheorem{definition}[theorem]{\sc Definition}
\newtheorem{remark}{\sc Remark}
\newtheorem{example}{\sc Example}
\newtheorem{definition-theorem}{{\bf{D}}{\fontsize{7pt}{11}\fontshape{n}
 \fontseries{bx}
 \selectfont{$\!\!\!$EFINITION~and~}}{\bf{T}}{\fontsize{7pt}{11}
 \fontshape{n}\fontseries{bx}\selectfont{$\!\!\!$HEOREM}}}[section]
\newtheorem{definition-proposition}[theorem]
 {{\bf{D}}{\fontsize{7pt}{11}\fontshape{n}
 \fontseries{bx}\selectfont{$\!\!\!$EFINITION~and~}}{\bf{P}}{\fontsize{7pt}{11}
 \fontshape{n}\fontseries{bx}\selectfont{$\!\!\!$ROPOSITION}}}
\newtheorem{definition-corollary}[theorem]
 {{\bf{D}}{\fontsize{7pt}{11}\fontshape{n}
 \fontseries{bx}\selectfont{$\!\!\!$EFINITION~and~}}{\bf{C}}{\fontsize{7pt}{11}
 \fontshape{n}\fontseries{bx}\selectfont{$\!\!\!$OROLLARY}}}
\newtheorem{definition-lemma}[theorem]
 {{\bf{D}}{\fontsize{7pt}{11}\fontshape{n}
 \fontseries{bx}\selectfont{$\!\!\!$EFINITION~and~}}{\bf{L}}{\fontsize{7pt}{11}
 \fontshape{n}\fontseries{bx}\selectfont{$\!\!\!$EMMA}}}
\newenvironment{proof}{\par\noindent{\sc Proof.}}{\endproof}
\renewcommand\endproof{\lfl$\sstyle\blacksquare$}

\newcommand\1{{1\mkern-5mu {\mathrm I}}}
\newcommand\abs{\par\vskip 0.3cm\goodbreak\noindent}

\newcommand\ad{\mathrm{ad}}

\newcommand\bbcbb{{{}^B_B{\mathcal C}^B_B}}
\newcommand\C{\mathcal{C}}

\newcommand\cpbi[1]{\mathop{{}_{\varphi_{1,2}}\!
                \overset{{}_{#1}}\bowtie_{\varphi_{2,1}}}}
\newcommand\cross[1]{\mathop{{}^{\nu_l}_{\mu_l}\!\overset{{}_{#1}}
            \bowtie\!{}^{\nu_r}_{\mu_r}}}
\newcommand\D{\mathcal{D}}

\newcommand\DY[1]{\mathcal{DY}\left({#1}\right)}

\newcommand\E{{1\mkern-5mu {\mathrm I}}}
\newcommand\End{\operatorname{End}}
\newcommand\field[1]{\mathbb{#1}}
\newcommand\hhchh{{{}^H_H{\mathcal C}^H_H}}

\newcommand\id{\mathrm{id}}

\newcommand\inj{\mathrm{i}}

\newcommand\lfl{\leaders\hbox to 1em{\hss \hss}\hfill}
\newcommand\m{\mathrm{m}}

\newcommand\nl{\par\noindent}
\newcommand\NN{\field{N}}

\newcommand\Obj{\operatorname{Obj}}
\newcommand\op{\mathrm{op}}
\newcommand\proj{\mathrm{p}}
\newcommand\QQ{\field{Q}}

\newcommand\scripthhchh{{{}^{\scriptscriptstyle H}_{\scriptscriptstyle H}
                {\mathcal C}^{\scriptscriptstyle H}_{\scriptscriptstyle H}}}

\newcounter{Part}
\begin{document}

\author{Yuri Bespalov \and Bernhard Drabant}
\title{\vskip -2truecm{\hfill\normalsize DAMTP-98-9}\\[2truecm]
  {\bfseries\LARGE Cross Product Bialgebras}\\
  {\bfseries\large Part \Roman{Part}}}
\date{February 1998}
\maketitle

\bibliographystyle{amsalpha}

\begin{abstract}
\noindent
The subject of this article are cross product bialgebras
without co-cycles. We establish a theory characterizing cross
product bialgebras universally in terms of projections and
injections. Especially all known types of biproduct, double cross 
product and bicross product bialgebras can be described by this 
theory. Furthermore the theory provides new families of
(co-cycle free) cross product bialgebras.
Besides the universal characterization we find an equivalent 
(co-)modular description of certain types of cross product bialgebras
in terms of so-called Hopf data. With the help of Hopf data construction
we recover again all known cross product bialgebras as well as new and
more general types of cross product bialgebras.
We are working in the
general setting of braided monoidal categories which allows us to
apply our results in particular to the braided category of Hopf bimodules
over a Hopf algebra. Majid's double biproduct is seen to be a 
twisting of a certain tensor product bialgebra in this category. This
resembles the case of the Drinfel'd double which can be
constructed as a twist of a specific cross product.
\abs
\end{abstract}

\section*{Introduction}

In recent years various (co-cycle free) cross products with bialgebra
structure had been investigated by several authors
\cite{Swe1:69,Rad1:85,Ma1:90}.
The different types like tensor product bialgebra, biproduct, 
double cross product and bicross product bialgebra are characterized each
by a universal formulation in terms of specific projections
and injections of the particular tensorands into the cross product.
The tensorands show interrelated (co-)module structures which are compatible
with these universal properties and which allow a reconstruction of the
cross product. The cross products are therefore equivalently
characterized by either of the two descriptions.
The multiplication and the comultiplication of the different cross products
have a similar form as the multiplication and the comultiplication of the
tensor product bialgebra except that the tensor transposition
is replaced by a more complicated morphism with
particular properties. The (co-)unit is given by the canonical
tensor product (co-)unit. Up to these common aspects the defining relations
of the several types of cross products seem to be 
different. The question arises if there exists at all a possibility
to describe the different cross products as different versions of a
single unifying theory which equivalently characterizes cross product
bialgebras universally and in a (co-)modular manner. 
The present article is concerned
with this question and will give an affirmative answer.

Based on the above mentioned common properties of cross products
we define cross product bialgebras or bialgebra
admissible tuples (BAT). We show that there are equivalent descriptions
of cross product bialgebras either by certain idempotents or by
coalgebra projections and
algebra injections obeying specific relations. However it is not clear
if a necessary and sufficient formulation by some interrelated
(co-)module structures of the particular tensor factors exists as well.
For these purposes we restrict the consideration to BATs where both the
algebra and coalgebra structure of the tensorands is respected at least by
one of the four projections or injections such that it is at the same
time an algebra and a coalgebra morphism. The corresponding cross
product bialgbras will be called
trivalent cross product bialgebras. This is a sufficiently general
class of cross product bialgebras to cover all the known cross products
of \cite{Swe1:69,Rad1:85,Ma1:90}. Trivalent cross product bialgebras
admit a universal characterization as well.
On the other hand we define so-called
Hopf data. A Hopf datum is a couple of objects which are both algebras and
coalgebras and which are mutual (co-)modules obeying certain
compatibility relations. One can show that a certain Hopf datum structure
is canonically inherited on any BAT. Conversely Hopf
data induce an algebra and a coalgebra on the tensor
product $B= B_1\otimes B_2$ which strongly resembles the
definition of a cross product bialgebra. However there is a priori no
compatibility of both structures rendering the Hopf datum a
bialgebra. However, Hopf data obey the fundamental recursive identities
$f=\Phi(f)$ of Proposition \ref{recurse-id} for both $f=\Delta_B\circ\m_B$ 
and $f=(\m_B\otimes\m_B)\circ(\id_B\otimes\Psi_{B,B}\otimes\id_B)\circ
(\Delta_B\otimes\Delta_B)$. This fact leads us to the definition of
so-called recursive Hopf data. Recursive Hopf data turn out to be
bialgebras. We also define recursive Hopf data with finite order and
show that a special class of recursive Hopf data (with order $\le 2$),
called trivalent Hopf data, are in one-to-one correpondence
with trivalent cross product bialgebras. Hence the classification of
cross product 
bialgebras either by (co-)modular or by universal properties according to
\cite{Rad1:85} has been achieved for all trivalent cross product
bialgebras in terms of trivalent Hopf data.
The new classification scheme covers all the known types
of cross products with bialgebra structure \cite{Swe1:69,Rad1:85,Ma1:90}.
And for the most general trivalent Hopf data it provides a new family of
cross product bialgebras which had not yet been studied in the literature
so far\footnote{A special variant of a trivalent Hopf datum has been
studied in \cite{Bes1:93} which uniformly describes biproducts and
bicross products in the symmetric category of vector spaces.}.
At the end of Section \ref{cross-prod} we will apply our results in
particular to Radford's 4-parameter Hopf algebra $H_{n,q,N,\nu}$
introduced in \cite{Rad1:94}. It turns out that it is a biproduct
bialgebra over the sub-Hopf group algebra $k C_N$ of $H_{n,q,N,\nu}$.

Since we are working throughout in very general types of braided
categories, we can apply our results to the special case of
the braided category of Hopf bimodules over a Hopf algebra
(possibly in a braided category, too). We demonstrate that Majid's
double biproduct \cite{Ma1:95} is a bialgebra twist of a certain tensor
product bialgebra in the category of Hopf bimodules. 
An example of double biproduct bialgebra is Lusztig's construction 
of the quantum group $\mathbf{U}$ \cite{Lus1:93}. 
 
A more thorough investigation of Hopf data and cross product bialgebras
in Hopf bimodule categories will be presented in a forthcoming work.
Another application of our results shows that the braided matched pair
formulation in terms of a certain pairing only works if the mutual
braiding of the two objects of the pair is involutive. This confirms in
some sense a similar observation in \cite{Ma1:95}.
\abs
In Section \ref{prelim} we give a survey of previous results, notations
and conventions which we need in the following. In particular we
recall outcomes of \cite{BD:95,BD:97}. We use graphical calculus
for braided categories.
The subsections on Hopf bimodules and twisting will be needed only in
Section \ref{applic}. Section \ref{cross-prod} is devoted to the main
subject of the article. We define bialgebra admissible tuples (BAT) or
cross product bialgebras, trivalent cross product bialgebras, and
(recursive) Hopf data (of finite order). It turns out that
Hopf data with a trivial (co-)action are
recursive and of order $\le 2$. They will be called trivalent Hopf data.
We show that trivalent Hopf data and trivalent cross product
bialgebras are equivalent. They generalize the known ``classical"
cross products which will be recovered as certain special examples.
The results of Section \ref{cross-prod} will be applied in Section
\ref{applic}. Using results of \cite{BD:95} we demonstrate that the double 
biproduct \cite{Ma1:95} can be
obtained as a bialgebra twist from the tensor product of two
Hopf bimodule bialgebras. We show that the braided version of the matched 
pairing \cite{Ma2:94} yields a matched pair if and only if the braiding
of the two tensorands is involutive. A matched pair is a special
kind of Hopf datum studied in Section \ref{cross-prod}.

\section{Preliminaries}\Label{prelim}

We presume reader's knowledge of the theory of braided monoidal categories.
Braided categories have been introduced in the work of Joyal
and Street \cite{JS:86,JS:93}. Since then they were studied intensively
by many authors. For an introduction to the theory of braided categories
we recommend to have a short look into the above mentioned articles or in
standard references on quantum groups and braided categories
\cite{CP:94,Tur1:94,Kas1:95,Ma4:93}.
Because of Mac Lane's Coherence Theorem for monoidal categories
\cite{Mac2:63,Mac1:72} we may restrict our conideration to strict braided
categories. In our article we denote
categories by caligraphic letters $\C$, $\D$, etc. For a braided monoidal
category $\C$ the tensor product is denoted by $\otimes_\C$, the unit
object by $\E_\C$, and the braiding by ${}^\C\Psi$. If it is clear from
the context we omit the index '$\C$' at the various symbols.
Henceforth we consider braided categories which admit split idempotents
\cite{BD:95,Lyu1:95,Lyu2:95}; for each idempotent $\Pi=\Pi^2:M\to M$ of any
object $M$ in $\C$ there exists an object $M_\Pi$ and a pair of morphisms
$(\inj_\Pi,\proj_\Pi)$ such that $\proj_\Pi\circ\inj_\Pi=\id_{M_\Pi}$
and $\inj_\Pi\circ\proj_\Pi=\Pi$. This is not a severe restriction
of the categories under consideration since every braided category
can be canonically embedded into a braided category which admits split
idempotents \cite{BD:95,Lyu1:95}.

We use and investigate algebraic structures in braided categories. We
suppose that the reader is familiar with the generalization of algebraic
structures to braided categories. Essentially we are working with
algebras, coalgebras, bialgebras, Hopf algebras, modules, comodules,
bimodules and bicomodules in braided categories
\cite{Lyu1:95,Ma4:93,Ma2:94}. Structures like Hopf
bimodules or crossed modules will be reviewed in the following.
We use throughout the symbol $\m$ for the multiplication and
$\eta$ for the unit of an algebra, $\Delta$ for the comultiplication and
$\varepsilon$ for the counit of a coalgebra, $S$ for the antipode of a
Hopf algebra, $\mu$ for the (left or right) action of an algebra on a
module, and $\nu$ for the (left or right) coaction of a coalgebra on a
comodule. We call a morphism $\rho:M\otimes N\to \E_\C$ in $\C$ a pairing
of $M$ and $N$. The graphical calculus for (strict) braided monoidal catgories
\cite{JS:86,FY:89,Kas1:95,ReT:90,ReT:91,Ma4:93,Tur1:94} will be used
throughout the paper.
We compose morphisms from up to down, i.~e.~the domains of the morphisms
are at the top and the codomaines are at the bottom of the graphics.
Tensor products are represented horizontally in the corresponding order.
We present our own conventions \cite{BD:95,BD:97} in Figure \ref{fig-conv}
and omit an assignment to a specific object if there is no fear of
confusion.
\begin{figure}\Label{fig-conv}
\begin{equation*}
\begin{array}{c}
\m=\divide\unitlens by 2
   \begin{tangle}\cu\end{tangle}
   \multiply\unitlens by 2
\quad ,\quad
\eta=\divide\unitlens by 2
   \ \begin{tangle}\unit\end{tangle}
   \multiply\unitlens by 2
\quad ,\quad
\Delta=\divide\unitlens by 2
   \begin{tangle}\cd\end{tangle}
   \multiply\unitlens by 2
\quad ,\quad
\varepsilon=\divide\unitlens by 2
   \ \begin{tangle}\counit\end{tangle}
   \multiply\unitlens by 2
\quad ,\quad
S=\divide\unitlens by 2
   \begin{tangle}\morph{S}\end{tangle}
   \multiply\unitlens by 2
\quad ,\quad
\\[8pt]
\mu_l=\divide\unitlens by 2
   \begin{tangle}\lu\end{tangle}
   \multiply\unitlens by 2
\quad ,\quad
\mu_r=\divide\unitlens by 2
   \begin{tangle}\ru\end{tangle}
   \multiply\unitlens by 2
\quad ,\quad
\nu_l=\divide\unitlens by 2
   \begin{tangle}\ld\end{tangle}
   \multiply\unitlens by 2
\quad ,\quad
\nu_r=\divide\unitlens by 2
   \begin{tangle}\rd\end{tangle}
   \multiply\unitlens by 2
\quad ,\quad
\\[8pt]
\Psi=\divide\unitlens by 2
   \begin{tangle}\x\end{tangle}
   \multiply\unitlens by 2
\quad ,\quad
\Psi^{-1}=\divide\unitlens by 2
   \begin{tangle}\xx\end{tangle}
   \multiply\unitlens by 2
\quad ,\quad\\[8pt]
\rho=\divide\unitlens by 2
   \begin{tangle}\Pairing{\sstyle\rho}\end{tangle}
   \multiply\unitlens by 2
\quad ,\quad
\rho^{\text{--}}=\divide\unitlens by 2
   \begin{tangle}\Pairing{\sstyle{\rho^{\text{--}}}}\end{tangle}
   \multiply\unitlens by 2
\quad .\quad
\end{array}
\end{equation*}
\caption{\small Graphical presentation of multiplication $\m$, unit
$\eta$, comultiplication $\Delta$, counit $\varepsilon$, antipode $S$,
left action $\mu_l$, right action $\mu_r$, left coaction $\nu_l$,
right coaction $\nu_r$, braiding $\Psi$, inverse braiding
$\Psi^{-1}$, pairing $\rho$ and convolution inverse
$\rho^{\text{--}}$ (if it exists).}
\end{figure}
To elucidate graphical calculus we will represent below
the bialgebra axiom of multiplicativity of the comultiplication
both in the ordinary way of composition and by graphical symbols.
\begin{equation*}
\begin{gathered}
\Delta\circ\m =(\m\otimes\m)\circ(\id\otimes\Psi\otimes\id)\circ
(\Delta\otimes\Delta)
\\[5pt]
\divide\unitlens by 5
\begin{tangle}
\cu\\
\cd
\end{tangle}
\,=\,
\begin{tangle}
\cd\step[2]\cd\\
\id\step[2]\x\step[2]\id\\
\cu\step[2]\cu
\end{tangle}
\end{gathered}
\end{equation*}

The results of the following subsections on Hopf bimodules and twisting
will be needed in Section \ref{applic}.
They are not relevant for the central part of the article presented in
Section \ref{cross-prod}.

\subsection{Hopf Bimodules}

Hopf bimodules over a bialgebra $B$ in $\C$ are
$B$-bimodules and $B$-bicomodules such that the
actions are bicomodule morphisms through the diagonal coactions
on tensor products of comodules and the canonical comodule
structure on $B$ \cite{BD:95}. $B$-Hopf bimodules 
and bimodule-bicomodule morphisms constitute the category
$\bbcbb$. For the symmetric category of $k$-vector spaces Hopf
bimodules have been introduced in \cite{Wor1:89} under the name
bicovariant bimodules.

Supose that $H$ is a Hopf algebra in $\C$. Then there exists
a tensor bifunctor rendering $\hhchh$ a (braided) monoidal category
\cite{BD:95,BD:97}. The proper formulation of the corresponding
theorem requires two auxiliary bifunctors $\odot$ and $\boxdot$ on
$\hhchh$. Two objects $X$ and $Y$ of $\hhchh$ yield the $H$-Hopf bimodule
$X\boxdot Y=X\otimes Y$ with diagonal left and right actions
$\mu_{d,l}^{X\boxdot Y}$ and $\mu_{d,r}^{X\boxdot Y}$, and with induced
left and right coactions $\nu_{i,l}^{X\boxdot Y}=\nu_l^X\otimes\id_Y$ and
$\nu_{i,r}^{X\boxdot Y}=\id_X\otimes\nu_r^Y$. The Hopf bimodule
$X\odot Y$ is obtained by categorical dualization of the previous structures.
For Hopf bimodule morphisms $f$ and $g$ we define $f\odot g=f\boxdot g=
f\otimes g$. Then the categories $(\hhchh, \odot)$ and
$(\hhchh, \boxdot)$ are semi-monoidal, i.~e.\ they are categories
which are almost monoidal, except that the unit object and the
relations involving it are not required. For the definition of
the braiding of $\hhchh$ we use the natural transformation
$\Theta:\displaystyle \odot\overset\bullet\to\boxdot{}^\op$ given by
$\Theta_{X,Y}:=(\mu_l^Y\otimes\mu_r^X)\circ
 (\id_H\otimes\Psi_{X,Y}\otimes\id_H)\circ
 (\nu_l^X\otimes\nu_r^Y)\,:\, X\otimes Y \to Y\otimes X$.
In the following theorem the braided monoidal structure of $\hhchh$ is
described \cite{BD:95,BD:97}.

\begin{theorem}
\Label{Hopf-br}
The category $\hhchh$ of Hopf bimodules over $H$ is monoidal.
The unit object is the canonical Hopf bimodule $H$, and the tensor
product $\otimes_H$ is uniquely defined (up to monoidal equivalence)
by one of the following equivalent conditions for any pair of $H$-Hopf
bimodules $X$ and $Y$.
\begin{itemize}
\item
The $H$-Hopf bimodule $X\otimes_H Y$ is the tensor product over $H$ of the
underlying modules, and the canonical morphism
$\lambda^H_{X,Y}:X\odot Y\to X\otimes_H Y\cong X\underset H\otimes Y$
is functorial in $\hhchh$, i.\ e.\
$\lambda^H:\odot\overset\bullet\to \underset H\otimes$.
\item
The $H$-Hopf bimodule $X\otimes_H Y$ is the cotensor product over $H$ of the
underlying comodules, and the canonical morphism
$\rho^H_{X,Y}:X\,\underset H\square \,Y\cong X\otimes_H Y\to
X\boxdot Y$ is functorial in $\hhchh$, i.\ e.\
$\rho^H:\underset H\square\overset\bullet\to \boxdot $.
\end{itemize}
The corresponding natural morphisms $\lambda^H$ and $\rho^H$ obey the
identity
\begin{equation}
\mathlabel{rho-lambda}
\rho^H_{X,Y}\circ\lambda^H_{X,Y}=
(\mu^{X}_{r}\otimes\mu^{Y}_{l})\circ
(\id_X\otimes\Psi_{H\,H}\otimes\id_Y)\circ
(\nu^X_r\otimes\nu^Y_l)\,.
\end{equation}
The category $\hhchh$ is pre-braided through the pre-braiding
${}^{{}^\scripthhchh}\!\Psi_{X,Y}$ uniquely defined by the condition
$\rho^H_{Y,X}\circ{}^\scripthhchh\Psi_{X,Y}\circ\lambda^H_{X,Y}=
\Theta_{X,Y}$. It is braided if the antipode of $H$ is an isomorphism
in $\C$.\endproof
\end{theorem}
\abs
Another concept closely related to Hopf bimodules are crossed modules
\cite{Wor1:89,Yet1:89,Bes1:97,BD:95}. The connection of both notions
had been studied in \cite{Wor1:89} and was reformulated in \cite{Sbg1:94}
for symmetric categories of modules over commutative rings. The general
investigation for braided categories which admit split idempotents can be
found in \cite{BD:95}.

A right crossed module over the Hopf
algebra $H$ is an object $M$ in $\C$ which is both right $H$-module and
right $H$-comodule such that the following identity holds.
\begin{equation}\mathlabel{crossed-mod}
\begin{array}{c}
\divide\unitlens by 3
\begin{tangle}
  \object{\sstyle M}\step[4]\object{\sstyle H}\\
  \rd\Step\cd\\
  \id\step\x\Step\id\\
  \ru\Step\cu\\
  \object{\sstyle M}\step[4]\object{\sstyle H}
 \end{tangle}
\multiply\unitlens by 3
\quad =\quad
\divide\unitlens by 3
\begin{tangle}
\object{\sstyle M}\step[2.5]\object{\sstyle H}\\
\id\Step\hcd\\
\x\step\id\\
\id\Step\k\\
\x\step\id\\
\id\Step\hcu\\
\object{\sstyle M}\step[2.5]\object{\sstyle H}
\end{tangle}
\multiply\unitlens by 3
\end{array}
\end{equation}
Identity \eqref{crossed-mod} is the graphical representation
of the equation
\begin{equation*}
\begin{split}
&(\mu_r\otimes\m_H)\circ(\id_M\otimes\Psi_{M,H}\otimes\id_H)\circ
(\nu_r\otimes\Delta_H)\\
&=(\id_M\otimes\m_H)\circ(\Psi_{H,M}\otimes\id_H)\circ
(\id_H\otimes\nu_r\circ\mu_r)\circ(\Psi_{M,H}\otimes\id_H)\circ
(\id_M\otimes\Delta_H)\,.
\end{split}
\end{equation*}
The right $H$-crossed modules and the corresponding module-comodule
morphisms form a category which is denoted by $\DY{\C}^H_H$.
In this notation $\mathcal{D}$ stands for Drinfel'd who introduced the
quantum double $D(H)$ of a Hopf algebra $H$, and $\mathcal{Y}$ stands
for Yetter who identified 
the category of representations of $D(H)$ with the category of 
$H$-crossed modules. Therefore crossed modules are sometimes called
Drinfel'd-Yetter modules or Yetter-Drinfel'd modules. 
$\DY{\C}^H_H$ is a monoidal category through the tensor product and the unit
object of $\C$, and the diagonal (co-)actions for tensor products
of crossed modules \cite{Bes1:97,BD:95}. 
In the following we will outline the relation of Hopf bimodules and
crossed modules \cite{BD:95,BD:97}. If $H$ is a Hopf algebra in the
category $\C$ and $(X,\mu_r,\mu_l,\nu_r,\nu_l)$ is an $H$-Hopf bimodule
then there exists an object ${}_HX$ such that
${}_HX\cong \E\underset{H}{\otimes} X$,
${}_HX\cong \E\underset{H}{\square} X$ and
${}_{\scriptscriptstyle X}\proj\circ{}_{\scriptscriptstyle X}\inj=
\id_{{}_HX}$ where ${}_{\scriptscriptstyle X}\proj:\E\otimes X\cong X\to
{}_HX\cong \E\underset{H}{\otimes}X$ and
${}_X\inj:\E\underset H\square X \cong {}_HX\to X \cong \E\otimes X$ are
the corresponding universal morphisms. The assignment
${}_H(\_\,):\hhchh\longrightarrow \DY{\C}^H_H$, which is given through
${}_H(X):= \left({}_HX,\,{}_{\scriptscriptstyle X}\proj\circ\mu_r\circ
({}_{\scriptscriptstyle X}\inj\otimes\id_H),\,
({}_{\scriptscriptstyle X}\proj\otimes\id_H)\circ\nu_r\circ
{}_{\scriptscriptstyle X}\inj\right)$
for an object $X$, and through ${}_H(f)= {}_{\scriptscriptstyle Y}\proj
\circ f\circ{}_{\scriptscriptstyle X}\inj$ for a Hopf bimodule morphism
$f:X\to Y$, defines a functor into the category of crossed modules.
Conversely a full inclusion functor $H\!\ltimes\!(\_\,):
{\DY{\C}}^H_H\to {}^H_H{\C}^H_H$ of the category of right $H$-crossed
modules into the category of $H$-Hopf bimodules is defined by
$H\ltimes (X) ={(H\otimes X,\mu_{i,l}^{H\otimes X}, \nu_{i,l}^{H\otimes X},
\mu_{d,r}^{H\otimes X}, \nu_{d,r}^{H\otimes X})}$
for any right crossed module $X$ and by $H\ltimes(f)= \id_H\otimes f$ for
any crossed module morphism $f$. The action $\mu_{i,l}^{H\otimes X}$ is
the left action induced by $H$ and $\mu_{d,r}^{H\otimes X}$ is the
diagonal action of the tensor product $H$-module. In the dual way the
coactions $\nu_{i,l}^{H\otimes X}$ and $\nu_{d,r}^{H\otimes X}$ are
defined. The following theorem holds \cite{BD:95}.

\begin{theorem}\Label{yd-hopfbi}
Let $H$ be a Hopf algebra in $\C$ with isomorphic antipode. Then the
categories $\DY{\C}^H_H$ and $\hhchh$ are
braided monoidal equivalent by
$\displaystyle \DY{\C}^H_H
\genfrac{}{}{0pt}{2}{\xrightarrow[\hphantom{{}_H(-)}]{H\ltimes\!(-)}}
                   {\xleftarrow[{}_H(-)]{\hphantom{H\ltimes\!(-)}}}
{}_H^H{\C}^H_H$.\endproof
\end{theorem}
\abs
If not otherwise mentioned we subsequently assume that the antipode
of the Hopf algebra $H$ is an isomorphism in $\C$.

\begin{remark}\Label{remark-pi}
{\normalfont A mirror symmetric result corresponding to
Theorem \ref{yd-hopfbi} holds for left $H$-crossed modules and
$H$-Hopf bimodules. Henceforth we will denote the idempotents
${}_{\scriptscriptstyle X}\inj\circ{}_{\scriptscriptstyle X}\proj$
and $\inj_{\scriptscriptstyle X}\circ\proj_{\scriptscriptstyle X}$ of
a Hopf bimodule $X$ by ${}_X\Pi$ and $\Pi_X$ respectively.
Explicitely it holds ${}_X\Pi=\mu_l^X\circ(S_H\otimes\id_X)\circ\nu_l^X$
and $\Pi_X=\mu_r^X\circ(\id_X\otimes S_H)\circ\nu_r^X$ \cite{BD:95}.}
\end{remark}
\abs
This observation leads to the following useful lemma.

\begin{lemma}\Label{aux-psi-inv}
Suppose that $X$ and $Y$ are $H$-Hopf bimodules and
$f,g:X\odot Y\to X\boxdot Y$ are Hopf bimodule morphisms.
Then the identity
\begin{equation}
\mathlabel{wor-braid1}
{}^{\hhchh}\Psi_{X,Y}\circ\lambda^H_{X,Y}\circ({}_X\Pi\otimes\Pi_Y)=
\lambda^H_{Y,X}\circ\Psi_{X,Y}\circ({}_X\Pi\otimes\Pi_Y)
\end{equation}
holds. The identity 
\begin{equation}\mathlabel{wor-braid2}
({}_X\Pi\otimes\Pi_Y)\circ f\circ ({}_X\Pi\otimes\Pi_Y)=
({}_X\Pi\otimes\Pi_Y)\circ g\circ ({}_X\Pi\otimes\Pi_Y)
\end{equation}
implies $f=g$.
\end{lemma}

\begin{proof}
The composition of both sides of \eqref{wor-braid1} with the monomorphism
$\rho_{Y,X}$ obviously leads to the identity
$\Theta_{Y,X}\circ({}_X\Pi\otimes\Pi_Y)
=(\mu^{X}_{r}\otimes\mu^{Y}_{l})\circ
(\id_X\otimes\Psi_{H\,H}\otimes\id_Y)\circ
(\nu^X_r\otimes\nu^Y_l)\circ({}_X\Pi\otimes\Pi_Y)$
which in turn coincides with $\Psi_{Y,X}\circ({}_X\Pi\otimes\Pi_Y)$.
To prove the second statement of the lemma observe that
any Hopf bimodule morphism $f:X\odot Y\to X\boxdot Y$ can be
expressed in terms of $f^\prime= f\circ ({}_X\Pi\otimes\Pi_Y)$ and
subsequently in terms of
$f^{\prime\prime}=({}_X\Pi\otimes\Pi_Y)\circ f\circ ({}_X\Pi\otimes\Pi_Y)$
in the following way
\begin{equation}\mathlabel{wor-braid4}
\hstretch 75 \vstretch 75
f=\enspace
\begin{tangle}
\step\object{X}\Step\object{Y}\\
\ld\Step\rd\\
\id\step{\hstr{100}\O{{}_X\Pi}}\Step{\hstr{100}\O{\Pi_Y}}\step\id\\
\lu\Step\ru\\
\step\tu{f}\\
\Step\object{X\otimes Y}
\end{tangle}
\enspace=\enspace
\begin{tangle}
\ld\Step\rd\\
\id\step\tu{f^\prime}\step\id\\
\Put(18,12)[rb]{\mu_{d,l}}\lu[2]\Put(2,12)[lb]{\mu_{d,r}}\ru[2]
\end{tangle}
\enspace=\enspace
\begin{tangle}
\Put(8,10)[rb]{\nu^X_l}\ld\step[4]\Put(2,10)[lb]{\nu^Y_r}\rd\\
\hh\id\step\Put(2,10)[lb]{\nu^X_{\ad,r}}\rd\Step%
 \Put(8,10)[rb]{\nu^Y_{\ad,l}}\ld\step\id\\
\id\step\id\step\x\step\id\step\id\\
\hh\id\step\x\Step\x\step\id\\
\id\step\id\step\tu{f^{\prime\prime}}\step\id\step\id\\
\hh\id\step\Put(18,2)[rb]{\mu_{i,l}}\lu[2]%
 \Put(2,2)[lb]{\mu_{i,r}}\ru[2]\step\id\\
\Put(28,8)[rt]{\mu_{d,l}}\lu[3]\Put(2,8)[lt]{\mu_{d,r}}\ru[3]
\end{tangle}
\end{equation}
where $\nu^X_{\ad,r}, \nu^Y_{\ad,l}$ are the (braided) adjoint coactions 
\cite{Bes1:97,BD:95}. The second identity of \eqref{wor-braid4} is
derived from the module properties of $f$. In a similar way one
obtains the third identity of \eqref{wor-braid4}.\end{proof}
\abs
Relation \eqref{wor-braid1} is the braided counterpart of Woronowicz's
definition of braiding of Hopf bimodules (see \cite{Wor1:89}).
\abs
Finally we recall the first part of the canonical transformation
procedure of bialgebras in $\hhchh$ into bialgebras in $\C$ \cite{BD:95}
which we need in the following.

\begin{proposition}\Label{hbb-bp}
Let $H$ be a Hopf algebra in $\C$. A bialgebra
$\underline B=(B,\underline\m_B,\underline\eta_B,
\underline\Delta_B,\underline\varepsilon_B)$ in ${}^H_H{\C}^H_H$ can be
turned into a bialgebra $B=(B,\m_B,\eta_B,\Delta_B,\varepsilon_B)$ in $\C$
where the structure morphisms are given by
\begin{equation}\mathlabel{bialg-trans2}
\m_B=\underline\m_B\circ\lambda^H_{B,B}\,,\quad
  \eta_B=\underline\eta_B\circ\eta_H\,,\quad
  \Delta_B=\rho^H_{B,B}\circ\underline\Delta_B\,,\quad
  \varepsilon_B=\varepsilon_H\circ\underline\varepsilon_B\,.
\end{equation}
If $\underline B=(B,\underline\m_B,\underline\eta_B,
\underline\Delta_B,\underline\varepsilon_B,\underline{S}_B)$ is Hopf
algebra in $\hhchh$ then $B=(B,\m_B,\eta_B,\Delta_B,\varepsilon_B,S_B)$
is Hopf algebra in $\C$ with antipode $S_B$ given by
$S_B=\underline S_B\circ S_{B/H}=S_{B/H}\circ \underline S_B$
where $S_{B/H}=\mu_l\circ(\id_H\otimes\mu_r)\circ(S_H\otimes\id_B\otimes
S_H)\circ(\id_H\otimes\nu_r)\circ\nu_l$.\endproof
\end{proposition}

\subsection*{Twisting}

In this subsection we present the twisting construction for bialgebras
in a braided category $\C$. We proceed along the lines of
\cite{Dri1:90,Ma6:94}.

Let $C$ be a coalgebra and $\chi:C\to\1_\C$ be a morphism into the unit
object. Henceforth we will use the following notations
\begin{equation}\mathlabel{conv-prod}
\chi.f:=(\chi\otimes f)\circ\Delta\,,\qquad
f.\chi:=(f\otimes\chi)\circ\Delta
\end{equation}
for any morphism $f:C\to B$ in $\C$.

\begin{definition}\Label{bialg-twist}
If $B$ is a bialgebra in $\C$ and $\chi:B\otimes B\to\1_\C$ is a morphism
obeying the identities
\begin{gather}
\mathlabel{2cocycle1}
\chi\circ(\id_B\otimes\chi.\m)=\chi\circ(\chi.\m\otimes\id_B)\,,
\\
\mathlabel{2cocycle2}
\chi\circ(\eta\otimes\id_B)=\varepsilon=\chi\circ(\id_B\otimes\eta)
\end{gather}
then $\chi$ is called a $2$-cocycle of the bialgebra $B$.
If $\chi$ is a convolution invertible $2$-cocycle, then the twist
$\m^\chi_B$ of the multiplication $\m_B$ is defined by
$\m^\chi_B:=\chi.\m_B.\chi^{\text{--}}$. If $B$ is a Hopf algebra then the
twist $S^\chi_B$ of the antipode $S$ is given by
$S^\chi_B:=u.S_B.u^{\text{--}}$ where
$u=\chi\circ(\id_B\otimes S_B)\circ\Delta_B$.
\end{definition}
\abs
\begin{remark} {\normalfont Under the condition of Definition
\ref{bialg-twist} the first identity in \eqref{2cocycle2} holds
if and only if the second one is valid.}
\end{remark}
\abs
In analogy to \cite{Dri1:90,Ma6:94} the following proposition can be
verified in the braided case because nowhere in the proof the involutivity
$\Psi^2=\id$ is needed. Therefore we will only sketch how to prove the
proposition.

\begin{definition-proposition}\Label{twist-bialg}
Let $B$ be bialgebra (Hopf algebra) and $\chi:B\otimes B\to\1_\C$ be an
invertible $2$-cocycle. Then $B_\chi:=
(B,\m^\chi_B,\eta_B,\Delta_B,\varepsilon_B,(S^\chi_B))$ with the twisted
multiplication $\m^\chi_B$ (and twisted antipode $S^\chi_B$) is again
bialgebra (Hopf algebra). $B_\chi$ is called the twisted bialgebra
(Hopf algebra) of $B$ obtained by the twist $\chi$.
\end{definition-proposition}

\begin{proof}
At first we will demonstrate that the bialgebra axiom
$L:=\Delta\circ\m=
(\m\otimes\m)\circ(\id\otimes\Psi\otimes\id)\circ(\Delta\otimes\Delta):=R$
for $B$ is equivalent to the bialgebra axiom $L_\chi=R_\chi$ for $B_\chi$.
This follows from the identities $(L_\chi).\chi=\chi.L$ and
$(R_\chi).\chi=\chi.R$. Secondly, using the previous fact we show that
the associativity $A^l:=\m\circ(\m\otimes\id)=\m\circ(\id\otimes\m)=:A^r$
of $B$ is equivalent to the associativity of $B_\chi$, denoted by
$A^l_\chi=A^r_\chi$. This is proved with the help of the identities
$(A^l_\chi).\big(\chi\circ(\id_A\otimes\chi.\m)\big)=
\big(\chi\circ(\chi.\m\otimes\id_A)\big).A^l$ and
$(A^r_\chi).\big(\chi\circ(\chi.\m\otimes\id_A)\big)=
\big(\chi\circ(\id_A\otimes\chi.\m)\big).A^r$ which result from
\eqref{2cocycle1}.
\end{proof}

\section{Cross Product Bialgebras and Hopf Data}\Label{cross-prod}

Section \ref{cross-prod} is the central part of the article. We define
cross product bialgebras or bialgebra admissible tuples (BAT) and Hopf
data. We consider certain specializations of these definitions, which we
call trivalent cross product bialgebras and recursive Hopf data
respectively. Trivalent cross product bialgebras form a sufficiently
general class to cover the cross product
bialgebras of \cite{Swe1:69,Rad1:85,Ma1:90}. Additionally there arise new
explicit examples of trivalent cross product bialgebras. All of them
will be completely classified in terms of recursive Hopf data.
Therefore an equivalent description either through interrelated
(co-)module structures or through universal projector decompositions is found.

\subsection{Cross Product Bialgebras} 

Suppose now there are two objects $B_1$ and $B_2$ in $\C$, and morphisms
$\varphi_{1,2}:B_1\otimes B_2\to B_2\otimes B_1$ and
$\varphi_{2,1}:B_2\otimes B_1\to B_1\otimes B_2$.

\begin{definition}\Label{bat}
We call $\big((B_1,\m_1,\eta_1,\Delta_1,\varepsilon_1),
(B_2,\m_2,\eta_2,\Delta_2,\varepsilon_2),\varphi_{1,2},\varphi_{2,1}\big)$
a bialgebra admissible tuple (BAT)
in the category $\C$ if $(B_j,\m_j,\eta_j)$ is an algebra and
$(B_j,\Delta_j,\varepsilon_j)$ is a coalgebra
for $j\in\{1,2\}$, such that $\varepsilon_j\circ\eta_j=\id_{\E_\C}$,
and the object $B_1\otimes B_2$
is a bialgebra through
\begin{equation}\mathlabel{bat-bialg}
\begin{split}
&\m_\times =(\m_1\otimes\m_2)\circ(\id_{B_1}\otimes\varphi_{2,1}\otimes
             \id_{B_2})\,,\quad
\eta_\times =\eta_1\otimes\eta_2\,,\\
&\Delta_\times = (\id_{B_1}\otimes\varphi_{1,2}\otimes\id_{B_2})\circ
            (\Delta_1\otimes\Delta_2)\,,\quad
\varepsilon_\times =\varepsilon_1\otimes\varepsilon_2\,.
\end{split}
\end{equation}
This bialgebra will be called the cross product bialgebra associated to the
bialgebra admissible tuple $(B_1,B_2,\varphi_{1,2},\varphi_{2,1})$ and is
denoted by $B_1\cpbi{} B_2$.
\end{definition}
\abs
One observes that the definition of a cross product bialgebra differs from the
usual definition of a canonical tensor product bialgebra (in a symmetric
category) only through the substitution of the tensor transposition by
the morphisms $\varphi_{1,2}$ and $\varphi_{2,1}$.

Of course the known cross products with bialgebra structure are cross
product bialgebras in the sense of Definition \ref{bat}
if the structures of the objects $B_1$ and $B_2$ and of the
morphisms $\varphi_{1,2}$ and $\varphi_{2,1}$ are chosen correctly.

The following proposition allows us to express cross product bialgebras
in terms of idempotents or projections and injections.

\begin{proposition}\Label{crossprod-proj}
Let $A$ be a bialgebra in $\C$, then the following statements are
equivalent.

\begin{enumerate}
\item
$A$ is bialgebra isomorphic to a cross product bialgebra $B_1\cpbi{} B_2$.
\item
There are idempotents $\Pi_1,\Pi_2\in\End(A)$ such that
\begin{equation*}
\begin{split}
&\m_A\circ(\Pi_j\otimes\Pi_j)=\Pi_j\circ\m_A\circ(\Pi_j\otimes\Pi_j)\,,
 \quad \Pi_j\circ\eta_A=\eta_A\,,\\
&(\Pi_j\otimes\Pi_j)\circ\Delta_A=(\Pi_j\otimes\Pi_j)\circ\Delta_A\circ
 \Pi_j\,,\quad \epsilon_A\circ\Pi_j=\epsilon_A
\end{split}
\end{equation*}
for $j\in\{1,2\}$, and the sequence
$A\otimes A\xrightarrow{\m_A\circ(\Pi_1\otimes\Pi_2)}A
         \xrightarrow{(\Pi_1\otimes\Pi_2)\circ\Delta_A} A\otimes A$
is a splitting of the idempotent $\Pi_1\otimes\Pi_2$ of $A\otimes A$.
\item
There exist objects $B_1$ and $B_2$ in $\C$ which are at the same time
algebras and coalgebras, and algebra morphisms $\inj_j$, coalgebra
morphisms $\proj_j$, $B_j\overset{\inj_j}\longrightarrow
A\overset{\proj_j}\longrightarrow B_j$, such that
$\proj_j\circ\inj_j=\id_{B_j}$ for $j\in\{1,2\}$, and the morphisms
$\m_A\circ(\inj_1\otimes \inj_2):B_1\otimes B_2\to A$ and
$(\proj_1\otimes \proj_2)\circ\Delta_A:A\to B_1\otimes B_2$
are inverse to each other.
\end{enumerate}
\end{proposition}

\begin{proof}
``(2) $\Rightarrow$ (3)'': Since $\Pi_j$ for $j\in\{1,2\}$ are
idempotents there are morphisms $\inj_j:B_j\to A$ and $\proj_j:A\to B_j$
which split $\Pi_j$. We define
\begin{equation*}
\begin{split}
\m_j &:= \proj_j\circ\m_A\circ(\inj_j\otimes\inj_j)\,,\\
\eta_j &:= \proj_j\circ\eta_A\,,\\
\Delta_j &:= (\proj_j\otimes\proj_j)\circ\Delta_A\circ\inj_j\,,\\
\varepsilon_j &:= \varepsilon_A\circ\inj_j
\end{split}
\end{equation*}
for $j\in\{1,2\}$. One immediately verifies that $(B_j,\m_j,\eta_j)$
are algebras and $(B_j,\Delta_j,\varepsilon_j)$ are coalgebras,
and $\inj_j$ are algebra morphisms, $\proj_j$ are coalgebra morphisms for
$j\in \{1,2\}$. Because $\Pi_j$ are idempotents it follows
$\m_A\circ(\inj_1\otimes\inj_2)\circ(\proj_1\otimes\proj_2)\circ\Delta_A
=\m_A\circ(\Pi_1\otimes\Pi_2)\circ(\Pi_1\otimes\Pi_2)\circ\Delta_A
=\id_A$ where the last equation follows by assumption. Similarly
$(\proj_1\otimes\proj_2)\circ\Delta_A\circ\m_A\circ(\inj_1\otimes\inj_2)=
\id_{B_1\otimes B_2}$ is proven.
\nl
``(3) $\Rightarrow$ (2)'': For $j\in\{1,2\}$ we consider the idempotents
$\Pi_j=\inj_j\circ\proj_j$. Statement (2) is then proven easily with the
help of the assumed properties of $\inj_j$ and $\proj_j$.
\nl
``(1) $\Rightarrow$ (3)'': Let $\phi:B_1\cpbi{} B_2\to A$ be the
isomorphism of bialgebras. Then in particular
\begin{equation}\mathlabel{varphi-unit}
\begin{split}
\varphi_{2,1}\circ(\eta_2\otimes\id_{B_1}) &= \id_{B_1}\otimes\eta_2\\
\varphi_{2,1}\circ(\id_{B_2}\otimes\eta_1) &= \eta_1\otimes\id_{B_2}
\end{split}
\end{equation}
and dually analogous for $\varphi_{1,2}$ and $\varepsilon_1$, $\varepsilon_2$.
We define the morphisms
\begin{equation}\mathlabel{ip-phi}
{\begin{split}
\inj_1 &:= \phi\circ(\id_{B_1}\otimes\eta_2)\,,\\
\inj_2 &:= \phi\circ(\eta_1\otimes\id_{B_2})\,,
\end{split}}
\qquad
{\begin{split}
\proj_1 &:= (\id_{B_1}\otimes\varepsilon_2)\circ\phi^{-1}\,,\\
\proj_2 &:= (\varepsilon_1\otimes\id_{B_2})\circ\phi^{-1}\,.
\end{split}}
\end{equation}
Using \eqref{ip-phi} one verifies without problems that $\inj_j$ are algebra
morphisms. In a dual manner it is proven that $p_j$ are coalgebra morphisms
for $j\in\{1,2\}$. Since $\varepsilon_1\circ\eta_1=\id_\E=
\varepsilon_2\circ\eta_2$ it follows $\proj_j\circ\inj_j=\id_{B_j}$,
$j\in\{1,2\}$.
Because $\phi$ is bialgebra isomorphism it holds
$\m_A\circ(\inj_1\otimes\inj_2) = \m_A\circ(\phi\otimes\phi)\circ
(\id_{B_1}\otimes\eta_2\otimes\eta_1\otimes\id_{B_2})= \phi\circ
\m_{B_1\cpbi{}B_2}\circ(\id_{B_1}\otimes\eta_2\otimes\eta_1\otimes\id_{B_2})
=\phi$ where \eqref{varphi-unit} has been used in the third equation.
Dually one obtains $(\proj_1\otimes\proj_2)\circ\Delta_A =\phi^{-1}$.
\nl
``(3) $\Rightarrow$ (1)'': By assumption the isomorphism
$\phi:=\m_A\circ(\inj_1\otimes\inj_2)$ induces a bialgebra structure on
$B:=B_1\otimes B_2$ through
\begin{equation}\mathlabel{phi-ind-bialg}
{\begin{split}
\m_B &= \phi^{-1}\circ\m_A\circ(\phi\otimes\phi)\,,\\
\eta_B &= \phi^{-1}\circ\eta_A\,,
\end{split}}
\qquad
{\begin{split}
\Delta_B &= (\phi^{-1}\otimes\phi^{-1})\circ\Delta_A\circ\phi\,,\\
\varepsilon_B &=\varepsilon_A\circ\phi\,.
\end{split}}
\end{equation}
We have to show that $B$ with the structure \eqref{phi-ind-bialg}
is a cross product bialgebra.
At first we prove that the (co-)units of $B$ are given by the tensor
products of the particular (co-)units of $B_1$ and $B_2$. It holds
$\eta_B=(\proj_1\otimes\proj_2)\circ\Delta_A\circ\eta_A=
 \proj_1\circ\eta_A\otimes\proj_2\circ\eta_A$ 
because $A$ is a bialgebra. Furthermore
$\eta_A=\inj_1\circ\eta_1=\inj_2\circ\eta_2$ is satisfied since 
$\inj_1$ and $\inj_2$ are algebra morphisms. Combining these two equations
then yields $\eta_B=\eta_1\otimes\eta_2$ since
$\proj_j\circ\inj_j =\id_{B_j}$ for $j\in\{1,2\}$ by assumption. Dually
$\varepsilon_B=\varepsilon_1\otimes\varepsilon_2$ can be proven.
Now we are going to prove that $B$ has the structure of a cross
product bialgebra if we use the morphisms
$\varphi_{1,2} = (\varepsilon_1\otimes\id_{B_2}\otimes
\id_{B_1}\otimes\varepsilon_2)\circ\Delta_B$ and
$\varphi_{2,1}=\m_B\circ(\eta_1\otimes\id_{B_2}\otimes\id_{B_1}\otimes
\eta_2)$. Thereto we need the auxiliary relations
\begin{equation}\mathlabel{auxil}
\begin{split}
(\id_{B_1}\otimes\m_2)\circ(\phi^{-1}\otimes\proj_2) &=
\phi^{-1}\circ\m_A\circ(\id_A\otimes\Pi_2)\,,\\
(\m_1\otimes\id_{B_1})\circ(\proj_1\otimes\phi^{-1}) &=
\phi^{-1}\circ\m_A\circ(\Pi_1\otimes\id_A)
\end{split}
\end{equation}
which can be proven easily because $\phi$ is a bialgebra isomorphism and
$\inj_1$ and $\inj_2$ are algebra morphisms. Using \eqref{auxil}
we show that $\m_B$ is indeed a multiplication of the form \eqref{bat-bialg}.
\begin{equation}\mathlabel{mult-bat}
\begin{split}
&(\m_1\otimes\m_2)\circ(\id_{B_1}\otimes\varphi_{2,1}\otimes\id_{B_2})\\
&=(\m_1\otimes\id_{B_2})\circ(\proj_1\otimes\phi^{-1})\circ
 (\id_A\otimes\m_A)\circ(\inj_1\otimes\inj_2\otimes\phi)\\
&=\phi^{-1}\circ\m_A\circ(\phi\otimes\phi)\\
&=\m_B\,.
\end{split}
\end{equation}
In the first equation of \eqref{mult-bat} two times \eqref{auxil} has been
used. In the second equality we applied $\proj_1\circ\inj_1=\id_{B_1}$ and
again \eqref{auxil}. The third equation of \eqref{mult-bat} holds by
definition. Similarly it can be shown by dualization that the
comultiplication is that of a cross product bialgebra, 
i.~e.~$\Delta_B= (\id_{B_1}\otimes\varphi_{1,2}\otimes\id_{B_2})\circ
(\Delta_1\otimes\Delta_2)$.
Hence $(B_1,B_2,\varphi_{1,2},\varphi_{2,1})$ is a BAT and
$B=B_1\cpbi{} B_2$ its corresponding cross product bialgebra.
\end{proof}
\abs
Since identities \eqref{varphi-unit} and their dual analogues for 
$\varphi_{1,2}$ and $\varepsilon_1$, $\varepsilon_2$ hold 
for cross product bialgebras, we immediately derive

\begin{corollary}\Label{alg-coalg-morph}
Let $B_1\cpbi{} B_2$ be a cross product bialgebra. Then 
$\eta_1\otimes\id_{B_2}$ and $\id_{B_1}\otimes\eta_2$ are algebra
morphisms, and $\varepsilon_1\otimes\id_{B_2}$, 
$\id_{B_1}\otimes\varepsilon_2$ are coalgebra morphisms.\endproof
\end{corollary}
\abs
It is not clear if the very general definition of a cross
product bialgebra is in one-to-one correspondence with a description
in terms of pairs of (co-)algebras with certain interrelated
compatible (co-)module structures. Hence a classification of
cross product bialgebras in the sense of \cite{Rad1:85} may not
succeed at this general level. But for reasons of classification and
reconstruction this aspect is important. The known
cross products with bialgebra structure \cite{Swe1:69,Rad1:85,Ma1:90}
admit such a description. For later use we will therefore define
trivalent cross product bialgebras as follows.

\begin{definition}\Label{tri-cross}
A cross product bialgebra $B_1\cpbi{} B_2$ is called trivalent
if at least one of the
morphisms $\eta_1\otimes\id_{B_2}$, $\id_{B_1}\otimes\eta_2$,
$\varepsilon_1\otimes\id_{B_2}$, $\id_{B_1}\otimes\varepsilon_2$ is both
an algebra and a coalgebra morphism. In a slight abuse of notation we
denote the corresponding bialgebra by $B_1\cpbi{3} B_2$ without
indication of the specific algebra-coalgebra morphism.
\end{definition}
\abs
In particular all cross products in \cite{Swe1:69,Rad1:85,Ma1:90} are
trivalent. Up to now we investigated universality of cross product
bialgebras. In the following subsection we study cross product
bialgebras from a (co-)modular point of view. 

\subsection{Hopf Data}

\begin{definition}\Label{hopf-pair}
A Hopf datum $(B_1,B_2,\mu_l,\nu_l,\mu_r,\nu_r)$ in $\C$ consists of
two objects $B_1$ and $B_2$ which are both counital algebras and unital
coalgebras, and $(B_1,\mu_l)$ is left $B_2$-module, $(B_1,\nu_l)$ is
left $B_2$-comodule, $(B_2,\mu_r)$ is right $B_1$-module, and
$(B_1,\nu_r)$ is right $B_1$-comodule obeying the identities
\begin{equation*}
\begin{array}{c}
\mu_r\circ(\eta_2\otimes\id_1)=\eta_2\circ\varepsilon_1
 =(\id_2\otimes\varepsilon_1)\circ\nu_l\,,
\quad
\varepsilon_2\circ\mu_r=
\varepsilon_2\otimes\varepsilon_1=\varepsilon_1\circ\mu_l\,,
\\
\mu_l\circ(\id_2\otimes\eta_1)=\eta_1\circ\varepsilon_2
 =(\varepsilon_2\otimes\id_1)\circ\nu_r\,,
\quad
\nu_r\circ\eta_2=\eta_2\otimes\eta_1=\nu_l\circ\eta_1\,,
\end{array}
\end{equation*}
\begin{equation*}
\begin{array}{c}
\divide\unitlens by 3
\begin{tangle}
\object{\sstyle B_1}\step[2.2]\object{\sstyle B_1}\\
\cu\\
\cd
\end{tangle}
\,\, =\,\,
 \begin{tangle}
\hstep\object{\sstyle B_1}\step[2.2]\object{\sstyle B_1}\\
      \hh \cd\step\cd \\
      \hh \id\hstep\hld\step\id\step\id \\
      \id\hstep\id\hstep\hx\step\id \\
      \hh \id\hstep\hlu\step\id\step\id \\
      \hh \cu\step\cu
\end{tangle}
\quad ,\quad
\begin{tangle}
\object{\sstyle B_2}\step[2.2]\object{\sstyle B_2}\\
\cu\\
\cd
\end{tangle}
\,\, =\,\,
\begin{tangle}
\hstep\object{\sstyle B_2}\step[2.2]\object{\sstyle B_2}\\
      \hh \cd\step\cd \\
      \hh \id\step\id\step\hrd\hstep\id \\
      \id\step\hx\hstep\id\hstep\id \\
      \hh \id\step\id\step\hru\hstep\id \\
      \hh \cu\step\cu
 \end{tangle}
\\
\\
\text{Algebra-coalgebra compatibility,}
\end{array}
\end{equation*}
\begin{equation*}
\begin{array}{c}
\divide\unitlens by 3
 \begin{tangle}
      \hh \cd\step\cd \\
      \id\step\hx\step\id \\
      \lu\step\ru \\
      \hh \ld\step\rd \\
      \id\step\hx\step\id \\
      \hh \cu\step\cu
 \end{tangle}
=
 \begin{tangle}
      \cd\step\cd \\
      \rd\step\hx\step\ld \\
      \id\step\hx\step\hx\step\id \\
      \ru\step\hx\step\lu \\
      \cu\step\cu
 \end{tangle}
\\
\\
\text{Module-comodule compatibility,}
\end{array}
\end{equation*}
\begin{equation*}
\begin{array}{c}
\divide\unitlens by 2
\begin{tangle}
\hh\hcu\hstep\id\\
\hstep\ru
\end{tangle}
\multiply\unitlens by 2
=
\divide\unitlens by 3
 \begin{tangle}
      \hh \id\step\cd\step\cd \\
      \id\step\id\step\hx\step\id \\
      \hh \id\step\lu\step\ru \\
      \Ru\dd \\
      \cu
 \end{tangle}
\multiply\unitlens by 3
\quad ,\quad
\divide\unitlens by 2
\begin{tangle}
\hh\id\hstep\hcu\\
\lu
\end{tangle}
\multiply\unitlens by 2
=
\divide\unitlens by 3
 \begin{tangle}
      \hh \cd\step\cd\step\id \\
      \id\step\hx\step\id\step\id \\
      \hh \lu\step\ru\step\id \\
      \step\d\Lu \\
      \Step\cu
 \end{tangle}
\\
\\
\text{Module-algebra compatibility,}
\end{array}
\end{equation*}
\begin{equation*}
\begin{array}{c}
\divide\unitlens by 2
\begin{tangle}
\hstep\rd\\
\hh\hcd\hstep\id
\end{tangle}
\multiply\unitlens by 2
=
\divide\unitlens by 3
 \begin{tangle}
      \cd \\
      \Rd\d \\
      \hh \id\step\ld\step\rd \\
      \id\step\id\step\hx\step\id \\
      \hh \id\step\cu\step\cu
 \end{tangle}
\multiply\unitlens by 3
\quad ,\quad
\divide\unitlens by 2
\begin{tangle}
\ld\\
\hh\id\hstep\hcd
\end{tangle}
\multiply\unitlens by 2
=
\divide\unitlens by 3
 \begin{tangle}
      \Step\cd \\
      \step\dd\Ld \\
      \hh \ld\step\rd\step\id \\
      \id\step\hx\step\id\step\id \\
      \hh \cu\step\cu\step\id
 \end{tangle}
\\
\\
\text{Comodule-coalgebra compatibility,}
\end{array}
\end{equation*}
\begin{equation*}
\begin{array}{c}
\divide\unitlens by 2
\begin{tangle}
\hstep\ru\\
\hh\hcd
\end{tangle}
\multiply\unitlens by 2
=
\divide\unitlens by 3
 \begin{tangle}
      \hcd\step\cd \\
      \id\step\hx\step\ld \\
      \ru\step\hx\step\id \\
      \cu\step\ru
 \end{tangle}
\multiply\unitlens by 3
\quad ,\quad
\divide\unitlens by 2
\begin{tangle}
\lu\\
\hh\hstep\hcd
\end{tangle}
\multiply\unitlens by 2
=
\divide\unitlens by 3
 \begin{tangle}
      \cd\step\hcd \\
      \rd\step\hx\step\id \\
      \id\step\hx\step\lu \\
      \lu\step\cu
 \end{tangle}
\\
\\
\text{Module-coalgebra compatibility,}
\end{array}
\end{equation*}
\begin{equation*}
\begin{array}{c}
\divide\unitlens by 2
\begin{tangle}
\hh\hcu\\
\hstep\rd
\end{tangle}
\multiply\unitlens by 2
=
\divide\unitlens by 3
 \begin{tangle}
      \cd\step\rd \\
      \rd\step\hx\step\id \\
      \id\step\hx\step\lu \\
      \hcu\step\cu
 \end{tangle}
\multiply\unitlens by 3
\quad ,\quad
\divide\unitlens by 2
\begin{tangle}
\hh\hstep\hcu\\
\ld
\end{tangle}
\multiply\unitlens by 2
=
\divide\unitlens by 3
 \begin{tangle}
      \ld\step\cd \\
      \id\step\hx\step\ld \\
      \ru\step\hx\step\id \\
      \cu\step\hcu
 \end{tangle}
\\
\\
\text{Comodule-algebra compatibility.}
\end{array}
\end{equation*}
\end{definition}
\abs
At first sight the defining relations of Hopf data seem to be rather
complicated and impenetrable. However all the compatibility identities
only relate the different (co-)algebra and (co-)module structures.
Besides there are two remarkable symmetries of the definition of Hopf data.
The first one is the usual categorical duality in conjunction with the
transformation ``$\m\leftrightarrow\Delta$",
``$\eta\leftrightarrow\varepsilon$", ``$\mu_l\leftrightarrow\nu_l$"
and ``$\mu_r\leftrightarrow\nu_r$". The second one is a kind of mirror
symmetry with respect to a vertical axis of the defining equations
considered as graphics in three dimensional space, followed by the
transformation of the indices ``$1\leftrightarrow 2$" and
``$l\leftrightarrow r$". This observation considerably simplifies subsequent
considerations and calculations. In a first step we recover canonical 
(co-)algebra structures of Hopf datum.

\begin{proposition}\Label{hp-alg-coalg}
Let $(B_1,B_2,\mu_l,\nu_l,\mu_r,\nu_r)$ be a Hopf datum. We define
\begin{equation}\mathlabel{hp-alg-coalg-id}
\begin{array}{c}
\phi_{1,2}=
\divide\unitlens by 3
 \begin{tangle}
      \ld\step\rd \\
      \id\step\hx\step\id \\
      \hh \cu\step\cu
 \end{tangle}
\multiply\unitlens by 3
\quad\text{and}\quad
\phi_{2,1}=
\divide\unitlens by 3
 \begin{tangle}
      \hh \cd\step\cd \\
      \id\step\hx\step\id \\
      \lu\step\ru \\
 \end{tangle}
\multiply\unitlens by 3
\end{array}
\end{equation}
Then $B= B_1\otimes B_2$ is both an algebra and a coalgebra
through the structure morphisms
\begin{equation}\mathlabel{hp-mult-comult}
\begin{split}
&\m_B =(\m_1\otimes\m_2)\circ(\id_{B_1}\otimes\phi_{2,1}\otimes
             \id_{B_2})\,,\quad
\eta_B =\eta_1\otimes\eta_2\,,\\
&\Delta_B = (\id_{B_1}\otimes\phi_{1,2}\otimes\id_{B_2})\circ
            (\Delta_1\otimes\Delta_2)\,,\quad
\varepsilon_B =\varepsilon_1\otimes\varepsilon_2\,.
\end{split}
\end{equation}
\end{proposition}

\begin{proof}
It is a well known fact that the necessary and sufficient conditions
for $(B,\m_B,\eta_B)$ being an algebra are given by the following equations.
\begin{equation}\mathlabel{cond-alg}
\begin{split}
\phi_{2,1}\circ(\m_2\otimes\id_{B_1}) &= (\id_{B_1}\otimes\m_2)\circ
 (\phi_{2,1}\otimes\id_{B_2})\circ(\id_{B_2}\otimes\phi_{2,1})\,\\
\phi_{2,1}\circ(\id_{B_2}\otimes\m_1) &= (\m_1\otimes\id_{B_2})\circ
 (\id_{B_1}\otimes\phi_{2,1})\circ(\phi_{2,1}\otimes\id_{B_1})\,\\
\phi_{2,1}\circ(\eta_2\otimes\id_{B_1}) &= \id_{B_1}\otimes\eta_2\,\\
\phi_{2,1}\circ(\id_{B_2}\otimes\eta_1) &= \eta_1\otimes\id_{B_2}\,.
\end{split}
\end{equation}
The verification of the third and fourth equation of \eqref{cond-alg} can
be done straightforwardly using \eqref{hp-alg-coalg-id} and the
defining relations of a Hopf datum. The second equation of \eqref{cond-alg}
will be proven graphically.
\begin{equation}\mathlabel{fig-proofassoc}
\begin{array}{c}
\vstretch 70 \hstretch 80
\begin{tangle}
\hstep\object{\sstyle B_2}\step\hstep
\object{\sstyle B_1}\step\object{\sstyle B_1}\\
\hh \hstep\id\step\hstep\cu \\
\hh \cd\step\cd \\
\id\step\hx\step\id \\
\lu\step\ru\\
\step\object{\sstyle B_1}\step\object{\sstyle B_2}
\end{tangle}
 =
\vstretch 70 \hstretch 50
\begin{tangle}
\hh \cd\hstep\cd\step\cd \\
\hh \id\step\id\hstep\id\hstep\hld\step\id\step\id \\
\id\step\id\hstep\id\hstep\id\hstep\hx\step\id \\
\hh \id\step\id\hstep\id\hstep\hlu\step\id\step\id \\
\hh \id\step\id\hstep\cu\step\cu \\
\id\step\hx\Step\id \\
\lu\step\Ru
\end{tangle}
 =
\vstretch 70 \hstretch 50
\begin{tangle}
\hh \step\hstep\cd\step\cd\step\cd \\
\hh \step\hstep\id\step\id\step\id\hstep\hld\step\id\step\id \\
\hstep\dd\step\hx\hstep\id\hstep\hx\step\id \\
\hh \cd\step\cd\hstep\id\hstep\hlu\step\id\step\id \\
\id\step\hx\step\id\hstep\hx\step\hddcu \\
\hh \lu\step\ru\dd\step\ru \\
\hh \step\id\step\lu\step\hstep\id \\
\step\cu\step\hstep\id
\end{tangle}
=
\vstretch 50 \hstretch 50
\begin{tangle}
\hcd\step\cd\Step\hstep\id \\
\id\step\hx\Step\d\step\hstep\id \\
\lu\step\hdcd\step\cd\hstep\id \\
\step\id\step\id\step\hx\step\ld\hstep\id \\
\step\id\step\ru\step\hx\step\id\hstep\id \\
\step\d\cu\step\ru\hcd \\
\Step\d\d\step\hx\step\id \\
\Step\step\d\lu\step\ru \\
\hh \Step\Step\cu\step\id
\end{tangle}
=
\vstretch 70 \hstretch 50
\begin{tangle}
\hh \cd\step\cd\step\id \\
\id\step\hx\step\id\step\id \\
\hh \lu\step\ru\step\id \\
\hh \step\id\hstep\cd\step\cd \\
\step\id\hstep\id\step\hx\step\id \\
\hh \step\id\hstep\lu\step\ru \\
\hh \step\id\step\dd\step\id \\
\hh \step\cu\step\hstep\id
\end{tangle}
\end{array}
\end{equation}
The first identity of \eqref{fig-proofassoc} uses the algebra-coalgebra
compatibility of Definition \ref{hopf-pair}. In the second equation
the module-algebra compatibility is used. The third equality holds
because $B_1$ is a left $B_2$-module, and the fourth identity is
true because of the module-coalgebra compatibility.
This proves the second equation of \eqref{cond-alg}. The first
identity of \eqref{cond-alg} can be verified similarly. Hence
$(B,\m_B,\eta_B)$ is an algebra. It will be proven dually that 
$(B,\Delta_B,\varepsilon_B)$ is coalgebra.\end{proof}
\abs
Under the conditions of Proposition \ref{hp-alg-coalg} one proves 
that $\eta_1\otimes\id_{B_2}$, $\id_{B_1}\otimes\eta_2$ are algebra 
morphisms, and $\varepsilon_1\otimes\id_{B_2}$, $\id_{B_1}\otimes\varepsilon_2$
are coalgebra morphisms. The result of Proposition \ref{hp-alg-coalg} 
strongly resembles the definiton of a BAT. However there is a priori
no compatibility of the algebra and the coalgebra structure of 
$B=B_1\otimes B_2$ rendering $B$ a cross product bialgebra. On the other 
hand the following proposition is easily proved.

\begin{proposition}\Label{bat-hp}
A BAT $(B_1,B_2,\varphi_{1,2},\varphi_{2,1})$ yields
a Hopf datum $(B_1,B_2,\mu_l,\nu_l,\mu_r,\nu_r)$ through
\begin{equation}\mathlabel{bat-act}
\begin{split}
\mu_l &=(\id_{B_1}\otimes\varepsilon_2)\circ\varphi_{2,1}\,,\quad
\mu_r=(\varepsilon_1\otimes\id_{B_2})\circ\varphi_{2,1}\\
\nu_l &=\varphi_{1,2}\circ(\id_{B_1}\otimes\eta_2)\,,\quad
\nu_r=\varphi_{1,2}\circ(\eta_1\otimes\id_{B_2})\,.
\end{split}
\end{equation}
Conversely according to eqs.~\eqref{hp-alg-coalg-id} in
Proposition \ref{hp-alg-coalg} the resulting Hopf
datum can be transformed into the cross product bialgebra $B_1\cpbi{}B_2$
because it holds $\phi_{1,2}=\varphi_{1,2}$ and $\phi_{2,1}=\varphi_{2,1}$.
\end{proposition}

\begin{proof}
Equations of the form \eqref{cond-alg} (and their dual form)
hold in particular for $\varphi_{2,1}$ (and $\varphi_{1,2}$). Hence the
(co-)module properties, the module-algebra and the comodule-coalgebra
compatibility of Definition \ref{hopf-pair} can be derived for the
(co-)actions in \eqref{bat-act}. The unital coalgebra structure and the
counital algebra structure of $B_1$ and $B_2$, as well as the relations in
Definition \ref{hopf-pair} involving the (co-)actions and the (co-)units
can be shown easily. By assumption $B_1\cpbi{} B_2$ is a bialgebra and
therefore it holds
\begin{equation}\mathlabel{bialg-ax-cp}
\Delta_\times\circ\m_\times =
(\m_\times\otimes\m_\times)\circ(\id_{B_1\otimes B_2}\otimes
\Psi_{B_1\otimes B_2,B_1\otimes B_2}\otimes\id_{B_1\otimes B_2})
\circ(\Delta_\times\otimes\Delta_\times)\,.
\end{equation}
Then one deduces the algebra-coalgebra compatibility of a Hopf datum by
either composing \eqref{bialg-ax-cp} with $\circ(\id_{B_1}\otimes
\eta_2\otimes\id_{B_1}\otimes\eta_2)$ and
$(\id_{B_1}\otimes\varepsilon_2\otimes\id_{B_1}\otimes\varepsilon_2)\circ$
or with $\circ(\eta_1\otimes\id_{B_2}\otimes\eta_1\otimes\id_{B_2})$ and
$(\varepsilon_1\otimes\id_{B_2}\otimes\varepsilon_1\otimes\id_{B_2})\circ\,$.
The module-comodule compatibility is derived from \eqref{bialg-ax-cp}
by composition with
$\circ(\eta_1\otimes\id_{B_2}\otimes\id_{B_1}\otimes\eta_2)$ and
$(\varepsilon_1\otimes\id_{B_2}\otimes\id_{B_2}\otimes\varepsilon_2)\circ\,$.
If one applies $\circ(\eta_1\otimes\id_{B_2}\otimes\eta_1\otimes\id_{B_2})$ and
$(\varepsilon_1\otimes\id_{B_2}\otimes\id_{B_1}\otimes\varepsilon_2)\circ$
to \eqref{bialg-ax-cp} one gets the first identity of the comodule-algebra
compatibility. The second equation of the comodule-algebra compatibility
is derived by composing \eqref{bialg-ax-cp} with
$\circ(\id_{B_1}\otimes\eta_2\otimes\id_{B_1}\otimes\eta_2)$ and
$(\varepsilon_1\otimes\id_{B_2}\otimes\id_{B_1}\otimes\varepsilon_2)\circ\,$.
The module-coalgebra compatibility is proven dually.
The composition of \eqref{bialg-ax-cp} with
$\circ(\id_{B_1}\otimes\eta_2\otimes\eta_1\otimes\id_{B_2})$ and
$(\varepsilon_1\otimes\id_{B_2}\otimes\id_{B_1}\otimes\varepsilon_2)\circ$
yields $\varphi_{1,2}=\phi_{1,2}$. The identity $\varphi_{2,1}=\phi_{2,1}$
is derived from \eqref{bialg-ax-cp} with
$\circ(\eta_1\otimes\id_{B_2}\otimes\id_{B_1}\otimes\eta_2)$ and
$(\id_{B_1}\otimes\varepsilon_2\otimes\varepsilon_1\otimes\id_{B_2})
\circ\,$. This concludes the proof of the proposition.
\end{proof}
\abs
Hence Hopf data are more general objects than BATs, and by Proposition
\ref{bat-hp} we may interpret cross product bialgebras as some subclass
of Hopf data. Subsequently we will show that two noteworthy recursive 
identities are satisfied for
Hopf data. Although these identities are rather involved, they possess the
above mentioned dual and mirror symmetries. To avoid complications we
represent the identities graphically.

\begin{definition}\Label{phi}
Let $(B_1,B_2,\mu_l,\nu_l,\mu_r,\nu_r)$ be a Hopf datum in $\C$ and
$f\in\End_\C(B_1\otimes B_2\otimes B_1\otimes B_2)$ be any endomorphism
in $\C$. Then the mapping
$\Phi:\End_\C(B_1\otimes B_2\otimes B_1\otimes B_2)\to
\End_\C(B_1\otimes B_2\otimes B_1\otimes B_2)$ is given by
\begin{equation}\label{recurse-element}
\Phi(f) \quad = \quad
\divide\unitlens by 3
\begin{tangle}
\hstep\object{\sstyle B_1}\step[3.5]\object{\sstyle B_2}%
\step[3.5]\object{\sstyle B_1}\step[3.5]\object{\sstyle B_2}\\
\hstep\id\step[2.5]\cd\step\cd\step[2.5]\id \\
\hstep\id\step[1.5]\cd\step\hx\step\cd\step[1.5]\id \\
\hstep\id\step[1.5]\rd\step\hx\step\hx\step\ld\step[1.5]\id \\
\hcd\step\id\step\hx\step\id\step\id\step\hx\step\id\step\hcd \\
\id\hstep\hld\step\lu\step\id\step\id\step\id\step\id
                                    \step\ru\step\hrd\hstep\id \\
\id\hstep\id\hstep\x\hstep\ffbox4{\sstyle f}\hstep\x\hstep\id
   \hstep\id \\
\id\hstep\hlu\step\ld\step\id\step\id\step\id\step\id
                                     \step\rd\step\hru\hstep\id \\
\hcu\step\id\step\hx\step\id\step\id\step\hx\step
                                      \id\step\hcu \\
\hstep\id\step[1.5]\ru\step\hx\step\hx\step\lu\step[1.5]\id \\
\hstep\id\step[1.5]\cu\step\hx\step\cu\step[1.5]\id \\
\hstep\id\step[2.5]\cu\step\cu\step[2.5]\id\\
\hstep\object{\sstyle B_1}\step[3.5]\object{\sstyle B_2}%
\step[3.5]\object{\sstyle B_1}\step[3.5]\object{\sstyle B_2}
\end{tangle}
\multiply\unitlens by 3
\end{equation}
\end{definition}
\abs
In the following proposition the fundamental recursive relation
for Hopf data will be derived.

\begin{proposition}\Label{recurse-id}
Let $(B_1,B_2,\mu_l,\nu_l,\mu_r,\nu_r)$ be a Hopf datum in $\C$ and
consider the endomorphisms $f_1:=\Delta_B\circ\m_B$ and
$f_2:=(\m_B\otimes\m_B)\circ(\id_B\otimes\Psi_{B,B}\otimes\id_B)\circ
(\Delta_B\otimes\Delta_B)$, where $\m_B$ and $\Delta_B$ are defined
according to \eqref{hp-mult-comult}. Then it holds
$f_1 = \Phi(f_1)$ and $f_2 = \Phi(f_2)$.
\end{proposition}

\begin{proof}
For $f_1=\Delta_B\circ\m_B$ we obtain the result through the following
(graphical) identities.
\begin{equation}
\mathlabel{recurse-id2}
\unitlens 8pt
\begin{tangle}
\object{\sstyle B_1}\step[1.5]\object{\sstyle B_2}\Step
\object{\sstyle B_1}\step[1.5]\object{\sstyle B_2}\\
\hh\id\step\cd\step\cd\step\id\\
\hh\id\step\id\step\x\step\id\step\id\\
\hh\hstr{200}\id\hstep\hlu\hstep\hru\hstep\id\\
\hh\hstr{200}\cucd\hstep\cucd\\
\hh\hstr{200}\id\hstep\hld\hstep\hrd\hstep\id\\
\hh\id\step\id\step\x\step\id\step\id\\
\hh\id\step\cu\step\cu\step\id\\
\object{\sstyle B_1}\step[1.5]\object{\sstyle B_2}\Step
\object{\sstyle B_1}\step[1.5]\object{\sstyle B_2}
\end{tangle}
\enspace=\enspace
\begin{tangle}
\hh\hstep\id\step\cd\Step\cd\step\id\\
\hstep\id\step\id\step\x\step\id\step\id\\
\hh\cd\hstep\lu\Step\ru\hstep\cd\\
\hh\id\hstep\hld\step\cd\step\cd\step\hrd\hstep\id\\
\hh\id\hstep\id\hstep\x\step\id\step\id\step\x\hstep\id\hstep\id\\
\hh\id\hstep\hlu\step\cu\step\cu\step\hru\hstep\id\\
\hh\cu\hstep\ld\Step\rd\hstep\cu\\
\hstep\id\step\id\step\x\step\id\step\id\\
\hh\hstep\id\step\cu\Step\cu\step\id
\end{tangle}
\quad =\quad
\unitlens 7pt
\begin{tangle}
\hstep\id\step[3.5]\cd\step\cd\step[3.5]\id\\
\hh\hstep\id\step[2.5]\sw1\Step\x\Step\nw1\step[2.5]\id\\
\hstep\id\step[1.5]\cd\step\hddcd\step\hdcd\step\cd\step[1.5]\id\\
\hh\hstep\id\step[1.5]\rd\step\x\step\id\step\id\step\x\step\ld
        \step[1.5]\id\\
\hh\cd\step\id\step\x\step\lu\step\ru\step\x\step\id\step\cd\\
\id\hstep\hld\step\lu\step\cu\step\cu\step\ru\step\hrd\hstep\id\\
\id\hstep\id\hstep\x\step[2]\id\step[3]\id\step[2]\x\hstep\id\hstep\id\\
\id\hstep\hlu\step\ld\step\cd\step\cd\step\rd\step\hru\hstep\id\\
\hh\cu\step\id\step\x\step\ld\step\rd\step\x\step\id\step\cu\\
\hh\hstep\id\step[1.5]\ru\step\x\step\id\step\id\step\x\step\lu
        \step[1.5]\id\\
\hstep\id\step[1.5]\cu\step\hdcu\step\hddcu\step\cu\step[1.5]\id\\
\hh\hstep\id\step[2.5]\nw1\Step\x\Step\sw1\step[2.5]\id\\
\hstep\id\step[3.5]\cu\step\cu\step[3.5]\id
\end{tangle}
\quad =\quad
\begin{tangle}
\hstep\id\step[2.5]\cd\step\cd\step[2.5]\id\\
\hstep\id\step[1.5]\cd\step\hx\step\cd\step[1.5]\id\\
\hh\hstep\id\step[1.5]\rd\step\x\step\x\step\ld\step[1.5]\id\\
\hh\cd\step\id\step\x\step\id\step\id\step\x\step\id\step\cd\\
\id\hstep\hld\step\lu\step\id\step\id\step\id\step\id%
                                    \step\ru\step\hrd\hstep\id\\
\id\hstep\id\hstep\x\hstep\dbox4{\Delta\!\circ\!\m}\hstep\x\hstep%
\id\hstep\id\\
\id\hstep\hlu\step\ld\step\id\step\id\step\id\step\id%
                                     \step\rd\step\hru\hstep\id\\
\hh\cu\step\id\step\x\step\id\step\id\step\x\step\id\step\cu\\
\hh\hstep\id\step[1.5]\ru\step\x\step\x\step\lu\step[1.5]\id\\
\hstep\id\step[1.5]\cu\step\hx\step\cu\step[1.5]\id\\
\hstep\id\step[2.5]\cu\step\cu\step[2.5]\id
\end{tangle}
\end{equation}
The first identity of \eqref{recurse-id2} uses the algebra-coalgebra
compatibilities of Definition \ref{hopf-pair}.
In the second identity of \eqref{recurse-id2}
we used the module-coalgebra and the comodule-algebra
compatibilities. (Co-)associativity is applied to derive the third
equation of \eqref{recurse-id2}.

For $f_2=(\m_B\otimes\m_B)\circ(\id_B\otimes\Psi_{B,B}\otimes\id_B)\circ
(\Delta_B\otimes\Delta_B)$ the proof is given by the following 
graphical equalities.
\begin{equation*}
\unitlens 7pt
\begin{tangle}
\step[.75]\object{\sstyle B_1}\step[2.5]\object{\sstyle B_2}\step[2.5]
\object{\sstyle B_1}\step[2.5]\object{\sstyle B_2}\\
\hh{\hstr{150}\cd}\step{\hstr{150}\cd}\step{\hstr{150}\cd}\step
        {\hstr{150}\cd}\\
\hh\id\hstep{\hstr{200}\hld\hstep\hrd}\hstep\id\step
        \id\hstep{\hstr{200}\hld\hstep\hrd}\hstep\id\\
\hh\id\hstep\id\step\x\step\id\hstep\id\step
        \id\hstep\id\step\x\step\id\hstep\id\\
\hh\id\hstep\cu\step\cu\hstep\x\hstep\cu\step\cu\hstep\id\\
\hh\id\step\id\Step\x\step\x\Step\id\step\id\\
\hh\id\hstep\cd\step\cd\hstep\x\hstep\cd\step\cd\hstep\id\\
\hh\id\hstep\id\step\x\step\id\hstep\id\step
        \id\hstep\id\step\x\step\id\hstep\id\\
\hh\id\hstep{\hstr{200}\hlu\hstep\hru}\hstep\id\step
        \id\hstep{\hstr{200}\hlu\hstep\hru}\hstep\id\\
\hh{\hstr{150}\cu}\step{\hstr{150}\cu}\step{\hstr{150}\cu}\step
        {\hstr{150}\cu}\\
\step[.75]\object{\sstyle B_1}\step[2.5]\object{\sstyle B_2}\step[2.5]
\object{\sstyle B_1}\step[2.5]\object{\sstyle B_2}
\end{tangle}
\enspace=\enspace
\begin{tangle}
\hh{\hstr{300}\cd}\step{\hstr{150}\cd}\step{\hstr{150}\cd}\step
        {\hstr{300}\cd}\\
\hh\id\step{\hstr{400}\hld}\step{\hstr{200}\hrd}\hstep\id\step\id
        \hstep{\hstr{200}\hld}\step{\hstr{400}\hrd}\step\id\\
\hh\id\step\id\Step\x\step\id\hstep\id\step\id\hstep\id\step\x\Step
        \id\step\id\\
\hh\id\hstep\cd\step\cd\hstep\cu\hstep\id\step\id\hstep\cu\hstep\cd
        \step\cd\hstep\id\\
\hh\id\hstep\id\step\id\step\hrd\hstep\id\step\id\step\x\step\id\step
        \id\hstep\hld\step\id\step\id\hstep\id\\
\hh\id\hstep\id\step\x\hstep\id\hstep\id\step\x\step\x\step\id\hstep
        \id\hstep\x\step\id\hstep\id\\
\hh\id\hstep\id\step\id\step\hru\hstep\id\step\id\step\x\step\id\step
        \id\hstep\hlu\step\id\step\id\hstep\id\\
\hh\id\hstep\cu\step\cu\hstep\cd\hstep\id\step\id\hstep\cd\hstep\cu
        \step\cu\hstep\id\\
\hh\id\step\id\Step\x\step\id\hstep\id\step\id\hstep\id\step\x\Step
        \id\step\id\\
\hh\id\step{\hstr{400}\hlu}\step{\hstr{200}\hru}\hstep\id\step\id
        \hstep{\hstr{200}\hlu}\step{\hstr{400}\hru}\step\id\\
\hh{\hstr{300}\cu}\step{\hstr{150}\cu}\step{\hstr{150}\cu}\step
        {\hstr{300}\cu}
\end{tangle}
\enspace=\enspace
\begin{tangle}
\hh\step\id\step[3]{\hstr{300}\cd}\step{\hstr{300}\cd}\step[3]\id\\
\hh{\hstr{200}\cd\hstep\cd}\Step\id\step\id\Step
        {\hstr{200}\cd\hstep\cd}\\
\hh\id\hstep{\hstr{300}\hld}\step{\hstr{300}\hrd}\hstep\d\step[1.5]\id
        \step\id\step[1.5]\dd\hstep{\hstr{300}\hld}\step
        {\hstr{300}\hrd}\hstep\id\\
\hh\id\hstep\id\hstep{\hstr{200}\hld}\step\id\hstep
        {\hstr{200}\hld\hstep\hrd}\hstep\id\step\id\hstep
        {\hstr{200}\hld\hstep\hrd}\hstep\id\step{\hstr{200}\hrd}
        \hstep\id\hstep\id\\
\hh\id\hstep\id\hstep\id\step\id\step\id\hstep\id\step\x\step\id
        \hstep\id\step\id\hstep\id\step\x\step\id\hstep\id\step\id
        \step\id\hstep\id\hstep\id\\
\hh\id\hstep\id\hstep\id\step\x\hstep\cu\step\cu\hstep\x\hstep\cu\step
        \cu\hstep\x\step\id\hstep\id\hstep\id\\
\hh\id\hstep\id\hstep\id\step\id\step\x\step\sw1\hstep\dd\step\d\hstep
        \nw1\step\x\step\id\step\id\hstep\id\hstep\id\\
\hh\id\hstep\id\hstep\id\step\hrd\hstep\id\step\cu\hstep\sw1\Step
        \nw1\hstep\cu\step\id\hstep\hld\step\id\hstep\id\hstep\id\\
\hh\id\hstep\id\hstep\x\hstep\id\hstep\id\step[1.5]\x\step[4]\x
        \step[1.5]\id\hstep\id\hstep\x\hstep\id\hstep\id\\
\hh\id\hstep\id\hstep\id\step\hru\hstep\id\step\cd\hstep\nw1\Step
        \sw1\hstep\cd\step\id\hstep\hlu\step\id\hstep\id\hstep\id\\
\hh\id\hstep\id\hstep\id\step\id\step\x\step\nw1\hstep\d\step\dd\hstep
        \sw1\step\x\step\id\step\id\hstep\id\hstep\id\\
\hh\id\hstep\id\hstep\id\step\x\hstep\cd\step\cd\hstep\x\hstep\cd\step
        \cd\hstep\x\step\id\hstep\id\hstep\id\\
\hh\id\hstep\id\hstep\id\step\id\step\id\hstep\id\step\x\step\id
        \hstep\id\step\id\hstep\id\step\x\step\id\hstep\id\step\id
        \step\id\hstep\id\hstep\id\\
\hh\id\hstep\id\hstep{\hstr{200}\hlu}\step\id\hstep
        {\hstr{200}\hlu\hstep\hru}\hstep\id\step\id\hstep
        {\hstr{200}\hlu\hstep\hru}\hstep\id\step{\hstr{200}\hru}
        \hstep\id\hstep\id\\
\hh\id\hstep{\hstr{300}\hlu}\step{\hstr{300}\hru}\hstep\dd\step[1.5]\id
        \step\id\step[1.5]\d\hstep{\hstr{300}\hlu}\step
        {\hstr{300}\hru}\hstep\id\\
\hh{\hstr{200}\cu\hstep\cu}\Step\id\step\id\Step
        {\hstr{200}\cu\hstep\cu}\\
\hh\step\id\step[3]{\hstr{300}\cu}\step{\hstr{300}\cu}\step[3]\id
\end{tangle}
\enspace=\enspace
\begin{tangle}
\hstep\id\step[2.5]\cd\step\cd\step[2.5]\id\\
\hstep\id\step[1.5]\cd\step\hx\step\cd\step[1.5]\id\\
\hh\hstep\id\step[1.5]\rd\step\x\step\x\step\ld\step[1.5]\id\\
\hh\cd\step\id\step\x\step\id\step\id\step\x\step\id\step\cd\\
\id\hstep\hld\step\lu\step\id\step\id\step\id\step\id%
                                    \step\ru\step\hrd\hstep\id\\
\id\hstep\id\hstep\x\hstep\dbox4{f_2}\hstep\x\hstep%
\id\hstep\id\\
\id\hstep\hlu\step\ld\step\id\step\id\step\id\step\id%
                                     \step\rd\step\hru\hstep\id\\
\hh\cu\step\id\step\x\step\id\step\id\step\x\step\id\step\cu\\
\hh\hstep\id\step[1.5]\ru\step\x\step\x\step\lu\step[1.5]\id\\
\hstep\id\step[1.5]\cu\step\hx\step\cu\step[1.5]\id\\
\hstep\id\step[2.5]\cu\step\cu\step[2.5]\id
\end{tangle}
\end{equation*}
where the first identity requires the algebra-coalgebra compatibility,
the second uses the module-algebra and the comodule-coalgebra
compatibility, as well as the (co-)module properties of $B_1$ and $B_2$.
In the third equation we applied associativity and again the (co-)module 
properties of $B_1$ and $B_2$.
\end{proof}
\abs
\begin{remark} {\normalfont Observe that we did not need the complete list of
defining relations of a Hopf datum for the deduction of Proposition
\ref{recurse-id}. We only needed that $B_1$ and $B_2$ are both algebras
and coalgebras, $B_1$ is a left $B_2$-(co-)module, $B_2$ is a right
$B_1$-(co-)module, the algebra-coalgebra compatibility, the
module-coalgebra compatibility, the comodule-algebra compatibility, the
module-algebra compatibility, and the comodule-coalgebra compatibility.
In particular we did not use the module-comodule compatibility.}
\end{remark}
\abs
From Propositions \ref{hp-alg-coalg} and \ref{bat} we conclude that
Hopf data are more general objects than BATs and therefore do not
correspond to them directly. In the following we will restrict our
considerations to so-called recursive Hopf data. They do not
necessarily imply universal characterization. For that reason a more special 
kind of recursive Hopf data, so-called trivalent Hopf data, will be 
introduced below.
We will verify that trivalent Hopf data and trivalent cross product
bialgebras are indeed equivalent notions which provide a 
(co-)modular and universal description of cross product bialgebras.
The resulting theory will turn out to be general enough to unify all
the ``classical" cross products of \cite{Swe1:69,Rad1:85,Ma1:90}. And
it will even generate a new family of cross product bialgebras which
can be described in this manner. The equivalent formulation
of cross product bialgebras either by interrelated (co-)module structures
or by certain universal projections and injections therefore will be
provided by our theory.

\begin{definition}\Label{recursive-hp}
Let $(B_1,B_2,\mu_l,\nu_l,\mu_r,\nu_r)$ be a Hopf datum in $\C$
and define the idempotent $\pi:=\eta_1\circ\varepsilon_1\otimes
\id_{B_2\otimes B_1}\otimes \eta_2\circ\varepsilon_2$.
Suppose that for every endomorphism $f$ of
$B_1\otimes B_2\otimes B_1\otimes B_2$ there exists a non-negative
integer $n\in \NN_0$ such that $\Phi^n(f)=\Phi^n(\pi\circ f\circ\pi)$.
Then the Hopf datum $(B_1,B_2,\mu_l,\nu_l,\mu_r,\nu_r)$ is
called recursive.
If $\Lambda:=\{n\in\NN_0\vert \Phi^n(f)=\Phi^n(\pi\circ f\circ\pi)\ \forall\,
f\in \End\C(B_1\otimes B_2\otimes B_1\otimes B_2)\}$ 
is a non-empty set such that $n_0:=\min\{n\in\Lambda\}$ exists, we call
$(B_1,B_2,\mu_l,\nu_l,\mu_r,\nu_r)$ a recursive Hopf datum of order $n_0$.
\end{definition}
\abs
A consequence of the structure of $\Phi$ is the following lemma.

\begin{lemma}\Label{1-order}
For every Hopf datum $(B_1,B_2,\mu_l,\nu_l,\mu_r,\nu_r)$ it holds
$\pi\circ\Phi(f)\circ\pi = \pi\circ f\circ \pi$ 
for any $f\in \End_\C(B_1\otimes B_2\otimes B_1\otimes B_2)$.
If $(B_1,B_2,\mu_l,\nu_l,\mu_r,\nu_r)$ is recursive of order $n$
then $\Phi^m(f)= \Phi^n(f)\quad\forall\,m\ge n$.
\end{lemma}

\begin{proof}
The first statement is verified straightforwardly using Hopf datum
properties and the structure of $\Phi$ according to Definition \ref{phi}.
Then it follows for a Hopf datum of order $n$ that
$\Phi^n(f)=\Phi^n(\pi\circ f\circ\pi)=\Phi^n(\pi\circ\Phi(f)\circ\pi)
= \Phi^n(\Phi(f))=\Phi^{n+1}(f)$.\end{proof}
\abs
A Hopf datum $(B_1,B_2,\mu_l,\nu_l,\mu_r,\nu_r)$ is \textit{non-trivial} if
$B_1$ or $B_2$ is not isomorphic to $\E_\C$.

\begin{lemma} A recursive Hopf datum $(B_1,B_2,\mu_l,\nu_l,\mu_r,\nu_r)$ 
of finite order is non-trivial if and only if its order is greater than 0.
\end{lemma}

\begin{proof}
If the order is 0 then $f=\pi\circ f\circ\pi$ and in particular for
$f=\id_{B_1\otimes B_2\otimes B_1\otimes B_2}$ one shows that
$\id_{B_i}=\eta_i\circ\varepsilon_i$. Therefore
$B_i\cong\E_\C$ for $i=\{1,2\}$. 
Conversely if $B_1$ and $B_2$ are isomorphic to $\E_\C$ then
$\pi=\id_{B_1\otimes B_2\otimes B_1\otimes B_2}$ and one concludes that
the order of the Hopf datum is 0.
\end{proof}

\begin{remark}
{\normalfont We will henceforth assume that Hopf data are non-trivial
so that the order is greater than 0 if it exists.}
\end{remark}
\abs
For special categories recursivity of a Hopf datum implies its
finite order.

\begin{proposition}
Suppose that $\C$ is a category of modules over a commutative ring
$\sigma$, and $(B_1,B_2,\mu_l,\nu_l,\mu_r,\nu_r)$ is a Hopf datum in $\C$
with $B_1$ and $B_2$ free $\sigma$-modules of finite rank.
Then $(B_1,B_2,\mu_l,\nu_l,\mu_r,\nu_r)$ is recursive if and only if
it is recursive of finite order.
\end{proposition}

\begin{proof}
Suppose that $(B_1,B_2,\mu_l,\nu_l,\mu_r,\nu_r)$ is recursive
and let $\{v_i\}_{i\in I}$ be a (finite) basis of
$B_1\otimes B_2\otimes B_1\otimes B_2$. Then every endomorphism
$f\in\End_\C(B_1\otimes B_2\otimes B_1\otimes B_2)$ can be written
as $f=\sum_{I,I} f_{i,j}\cdot\hat\delta_{i,j}$ where
$f(v_i)=\sum_I f_{i,j}\cdot v_j$ and $\hat\delta_{i,j}(v_k)=\delta_{i,k}
\cdot v_j$. Define $n_0:=\min\{n\in\NN_0\vert \Phi^n(\hat\delta_{i,j})=
\Phi^n(\pi\circ\hat\delta_{i,j}\circ\pi)\ \forall\,i,j\in I\}$ which
exists since $\vert I\vert$ is finite and the Hopf datum is supposed to
be recursive. Then one immediately concludes that the order of the Hopf
datum is $n_0$.\end{proof}
\abs
The most striking aspect of recursive Hopf data is the bialgebra
structure of the corresponding tensor product (co-)algebra
$B_1\otimes B_2$.

\begin{theorem}\Label{recurse-hp-bialg}
If $(B_1,B_2,\mu_l,\nu_l,\mu_r,\nu_r)$ is a recursive Hopf datum then
$B=B_1\otimes B_2$ equipped with the structure morphisms $\m_B$, $\eta_B$,
$\Delta_B$ and $\varepsilon_B$ according to Proposition \ref{hp-alg-coalg}
is a bialgebra. It will be denoted by $B=B_1\cross{} B_2$.
\end{theorem}

\begin{proof}
Because of Proposition \ref{hp-alg-coalg} we only have to prove
that $\Delta_B$ is an algebra morphism. With the help of Proposition 
\ref{recurse-id} and the recursivity of the Hopf datum we derive
for some $n\in\NN$ the identities 
$f_1=\Phi^n(f_1)=\Phi^n(\pi\circ f_1\circ\pi)$ and
$f_2=\Phi^n(f_2)=\Phi^n(\pi\circ f_2\circ\pi)$
where $f_1=\Delta_B\circ\m_B$ and $f_2=(\m_B\otimes\m_B)\circ
(\id_B\otimes\Psi_{B,B}\otimes\id_B)\circ(\Delta_B\otimes\Delta_B)$.
The relation $\pi\circ f_1\circ\pi=\pi\circ f_2\circ\pi$
holds by the module-comodule compatibility of the Hopf data.
Therefore $f_1=f_2$ which proves the theorem.
\end{proof}
\abs
In the next proposition we will introduce trivalent Hopf data.

\begin{definition-proposition}\Label{trivalent-hp-bialg}
Suppose that $(B_1,B_2,\mu_l,\nu_l,\mu_r,\nu_r)$ is a Hopf datum in $\C$.
Then one of the (co-)actions of $(B_1,B_2,\mu_l,\nu_l,\mu_r,\nu_r)$
is trivial if and only if one of the morphisms $\eta_1\otimes\id_{B_2}$,
$\id_{B_1}\otimes\eta_2$, $\varepsilon_1\otimes\id_{B_2}$ or
$\,\id_{B_1}\otimes\varepsilon_2$ is both algebra and coalgebra
morphism. If these equivalent conditions hold, we call
$(B_1,B_2,\mu_l,\nu_l,\mu_r,\nu_r)$ a trivalent Hopf datum.
Trivalent Hopf data are recursive of order $\le 2$. Therefore the
corresponding bialgebra $B_1\cross{}B_2$ exists and will be denoted by
$B_1\cross{3}B_2$.
\end{definition-proposition}

\begin{proof}
Suppose that $\mu_l$ is trivial. Then from \eqref{hp-mult-comult} we
conclude $\m_B=(\m_1\otimes\m_2)\circ\big(\id_{B_1}\otimes
(\id_{B_1}\otimes\mu_r)\circ(\Psi_{B_2,B_1}\otimes\id_{B_1})\circ
(\id_{B_2}\otimes\Delta_1)\otimes\id_{B_2}\big)$. Therefore
$(\id_{B_1}\otimes\varepsilon_2)\circ\m_B =\m_1\circ
 (\id_{B_1}\otimes\varepsilon_2\otimes\id_{B_1}\otimes\varepsilon_2)$
and $(\id_{B_1}\otimes\varepsilon_2)\circ\eta_B = \eta_1$
which shows that $(\id_{B_1}\otimes\varepsilon_2)$ is an algebra morphism.
Because of Proposition \ref{hp-alg-coalg} $(\id_{B_1}\otimes\varepsilon_2)$
is also a coalgebra morphism.
Conversely if $(\id_{B_1}\otimes\varepsilon_2)$ is an algebra morphism
it holds in particular $(\id_{B_1}\otimes\varepsilon_2)\circ\m_B\circ
(\eta_1\otimes\id_{B_2}\otimes\id_{B_1}\otimes\eta_2)
=\m_1\circ(\id_{B_1}\otimes\varepsilon_2\otimes
\id_{B_1}\otimes\varepsilon_2)\circ
(\eta_1\otimes\id_{B_2}\otimes\id_{B_1}\otimes\eta_2)$
from which the triviality $\mu_l=\varepsilon_2\otimes\id_{B_1}$ of
$\mu_l$ is derived. Since $\mu_l$ is trivial the identity
\begin{equation}\mathlabel{recurse-hp1}
\Phi(f)=\Phi\big((\id_{B_1\otimes B_2\otimes B_1}\otimes\eta_2\circ
\varepsilon_2)\circ f\circ(\eta_1\circ\varepsilon_1\otimes\id_{B_2\otimes
B_1\otimes B_2})\big)
\end{equation}
can be verified directly for any endomorphism $f$. From the general structure
of $\Phi$ as given in Definition \ref{phi} we derive
\begin{align}
\mathlabel{recurse-hp2}
&(\id_{B_1\otimes B_2\otimes B_1}\otimes\eta_2\circ\varepsilon_2)\circ
\Phi(f)\circ(\eta_1\circ\varepsilon_1\otimes
\id_{B_2\otimes B_1\otimes B_2})\\
=\,&(\id_{B_1\otimes B_2\otimes B_1}\otimes\eta_2\circ\varepsilon_2)\circ
\Phi\big((\eta_1\circ\varepsilon_1\otimes\id_{B_2\otimes B_1\otimes B_2})
\circ f\circ(\id_{B_1\otimes B_2\otimes B_1}\otimes\eta_2\circ\varepsilon_2)
\big)\circ\notag\\
&\circ(\eta_1\circ\varepsilon_1\otimes
\id_{B_2\otimes B_1\otimes B_2})\,.\notag
\end{align}
Using \eqref{recurse-hp1} two times, \eqref{recurse-hp2} and
then again \eqref{recurse-hp1} one obtains the result
$\Phi^2(f) = \Phi^2(\pi\circ f\circ\pi)$
for any $f\in\End_\C(B_1\otimes B_2\otimes B_1\otimes B_2)$. Hence the
order of the Hopf datum is $\le 2$.
Because of the dual and mirror symmetries of Hopf data the proof
of the proposition follows analogously for any other (co-)action being 
trivial.\end{proof}
\abs
\begin{remark}\Label{non-triv-coact}
{\normalfont One verifies easily that \eqref{recurse-hp1}
can be obtained under the following conditions which are weaker than the
assumption of triviality of one of the (co-)actions.
\begin{equation*}
{\hstretch 50 \vstretch 50
\begin{tangle}
\object{\sstyle B_2}\step[2]\object{\sstyle B_1}\\
\rd\step\id\\
\hh\id\step\x\\
\lu\step\id
\end{tangle}}
\enspace=\enspace
{\hstretch 50 \vstretch 50
\begin{tangle}
\object{\sstyle B_2}\step[1.5]\object{\sstyle B_1}\\
\step[0.2]\lu\step\unit
\end{tangle}}
\qquad\text{and}\qquad
{\hstretch 50 \vstretch 50
\begin{tangle}
\step[-0.2]\object{\sstyle B_2}\step[1.5]\object{\sstyle B_2}\\
\id\step\rd\\
\hh\x\step\id\\
\id\step\lu
\end{tangle}}
\enspace=\enspace
{\hstretch 50 \vstretch 50
\begin{tangle}
\object{\sstyle B_2}\step[1.5]\object{\sstyle B_2}\\
\counit\step[1.5]\rd
\end{tangle}}
\end{equation*}
or similar conditions involving $\nu_l$ and $\mu_r$. Then it follows
quite analogously as in Proposition \ref{trivalent-hp-bialg} that the Hopf
datum is recursive of order $\le 2$.
Using the notion of Proposition \ref{hp-alg-coalg} these conditions
can be reformulated as conditions of $\phi_{1,2}$ and
$\phi_{2,1}$. These are conditions of a BAT since the Hopf datum is
recursive.  Therefore a generalization of Theorem \ref{rhp-cpb} below can
be obtained in this way (see also the notion ``strong Hopf datum'' in
\cite{BD:98}).} 
\end{remark}
\abs
The following theorem demonstrates that
trivalent Hopf data and trivalent cross product bialgebras coincide. This
shows that a unified theory of cross product bialgebras has
been found which provides universal and (co-)modular
characterization equivalently.

\begin{theorem}\Label{rhp-cpb}
Suppose that $A$ is a bialgebra in $\C$. Then the following statements
are equivalent.

\begin{enumerate}
\item
There is a trivalent Hopf datum $(B_1,B_2,\mu_l,\nu_l,\mu_r,\nu_r)$ such
that the corresponding bialgebra $B_1\cross{3} B_2$ is bialgebra isomorphic 
to $A$.
\item
$A$ is bialgebra isomorphic to a trivalent cross product bialgebra 
$B_1\cpbi{3} B_2$.
\item
There are algebra morphisms $\inj_j:B_j\to A$ and
coalgebra morphisms $\proj_j:A\to B_j$ such that
$\proj_j\circ\inj_j=\id_{B_j}$ for $j\in\{1,2\}$,
$\m_A\circ(\inj_1\otimes\inj_2)=
\big((\proj_1\otimes\proj_2)\circ\Delta_A\big)^{-1}$, and one of the
morphisms $\inj_1$, $\inj_2$, $\proj_1$, $\proj_2$ is both
algebra and coalgebra morphism.
\item
There are idempotents $\Pi_1,\Pi_2\in\End(A)$ such that
\begin{equation*}
\begin{split}
&\m_A\circ(\Pi_j\otimes\Pi_j)=\Pi_j\circ\m_A\circ(\Pi_j\otimes\Pi_j)\,,
 \quad \Pi_j\circ\eta_A=\eta_A\,,\\
&(\Pi_j\otimes\Pi_j)\circ\Delta_A=(\Pi_j\otimes\Pi_j)\circ\Delta_A\circ
 \Pi_j\,,\quad \epsilon_A\circ\Pi_j=\epsilon_A
\end{split}
\end{equation*}
for every $j\in\{1,2\}$, the sequence
$A\otimes A\xrightarrow{\m_A\circ(\Pi_1\otimes\Pi_2)}A
         \xrightarrow{(\Pi_1\otimes\Pi_2)\circ\Delta_A} A\otimes A$
is a splitting of the idempotent $\Pi_1\otimes\Pi_2$ of $A\otimes A$, and one
of the idempotents $\Pi_1$, $\Pi_2$ is either algebra or coalgebra morphism.
\end{enumerate}
\end{theorem}

\begin{proof}
Essentially the proof of the theorem has been done in Propositions
\ref{crossprod-proj}, \ref{hp-alg-coalg}, \ref{bat-hp}, and
\ref{trivalent-hp-bialg}. We only have to show the additional (co-)algebra
properties of the corresponding morphisms or the triviality of the
corresponding (co-)actions.
\nl
``(3) $\Rightarrow$ (1)": From Proposition \ref{crossprod-proj} it follows
especially that $A$ is isomorphic to a cross product bialgebra
$B=B_1\cpbi{} B_2$ through the bialgebra isomorphism
$\phi:B_1\cpbi{} B_2\to A$, $\phi=\m_A\circ(\inj_1\otimes\inj_2)$.
Then there is a Hopf datum such that $B=B_1\cross{} B_2=B_1\cpbi{} B_2$
because of Proposition \ref{bat-hp}.
Suppose that $\proj_1$ is algebra and coalgebra morphism. Then
$\proj_1\circ\phi\circ\m_B = \m_1\circ(\proj_1\otimes\proj_1)\circ
    (\inj_1\otimes\inj_2)\circ\m_B
  = \proj_1\circ\m_A\circ(\phi\otimes\phi)$.
Therefore $\m_1\circ(\proj_1\circ\inj_1\otimes\proj_1\circ\inj_2)\circ\m_B
=\m_1\circ(\m_1\otimes\m_1)\circ(\proj_1\circ\inj_1\otimes
\proj_1\circ\inj_2\otimes\proj_1\circ\inj_1\otimes\proj_1\circ\inj_2)$.
Since $\proj_1\circ\inj_1=\id_{B_1}$ and $\proj_1\circ\inj_2=
(\id_{B_1}\otimes\varepsilon_2)\circ\phi^{-1}\circ\phi\circ
(\eta_1\otimes\id_{B_2})=\eta_1\circ\varepsilon_2$ the
identity $\mu_l=(\id_{B_1}\otimes\varepsilon_2)\circ\phi=\varepsilon_2
\otimes\id_{B_1}$ follows then from \eqref{hp-mult-comult}
and \eqref{bat-act}.
\nl
``(1) $\Rightarrow$ (2)": If there is a trivalent Hopf datum with
trivial left action $\mu_l$ then $\mu_l=(\id_{B_1}\otimes\varepsilon_2)
\circ\phi=\varepsilon_2\otimes\id_{B_1}$ and therefore
$(\id_{B_1}\otimes\varepsilon_2)\circ\m_B =\m_1\circ(\id_{B_1}\otimes
\varepsilon_2\otimes\id_{B_1}\otimes\varepsilon_2)$ and
$(\id_{B_1}\otimes\varepsilon_2)\circ\eta_B =\eta_1$
which implies that $(\id_{B_1}\otimes\varepsilon_2)$ is an algebra
morphism. Because of Proposition \ref{hp-alg-coalg} it is also a
coalgebra morphism.
\nl
``(2) $\Rightarrow$ (4)": Suppose that $(\id_{B_1}\otimes\varepsilon_2)$
is an algebra morphism. Using Proposition \ref{crossprod-proj} and
the bialgebra isomorphism $\phi: B_1\cpbi{}B_2\to A$ we obtain
$\Pi_1 =\inj_1\circ\proj_1 = \phi\circ(\id_{B_1}\otimes\eta_2)\circ
 (\id_{B_1}\otimes\varepsilon_2)\circ\phi^{-1}$
and therefore $\Pi_1$ is an algebra morphism.
\nl
``(4) $\Rightarrow$ (3)": If $\Pi_1$ is an algebra morphism then
$\Pi_1\circ\m_1=\inj_1\circ\proj_1\circ\m_A=\m_A\circ(\inj_1\circ\proj_1
\otimes\inj_1\circ\proj_1)=\inj_1\circ\m_1\circ(\proj_1\otimes\proj_1)$
and $\Pi_1\circ\eta_A=\inj_1\circ\proj_1\circ\eta_A=\eta_A=\inj_1\circ\eta_1$
because $\inj_1$ is an algebra morphism. Since $\inj_1$ is monomorphic
one concludes that $\proj_1$ is an algebra morphism and
also a coalgebra morphism by Proposition \ref{crossprod-proj}.
Thus the Theorem is proven for a particular case. Because of dual and
mirror symmetry all other cases can be verified analogously.
\end{proof}
\abs
Theorem \ref{rhp-cpb} shows that there exists a one-to-one
correspondence of trivalent cross product bialgebras and trivalent Hopf
data.
In addition both notations are equivalent to a description in terms of
a certain projector decomposition. Since Definition \ref{tri-cross} of a
trivalent cross product bialgebra is a generalization of the cross
products with bialgebra structure according to
\cite{Swe1:69,Rad1:85,Ma1:90}, we can express all of them in a unified
manner through trivalent Hopf data. Moreover the most general trivalent
Hopf data give rise to a new family of (trivalent) cross product
bialgebras.

For a better understanding of the theory we list five special
examples of trivalent Hopf data in the sequel which cover all the other
cases because of the dual and mirror
symmetries. The discussion of the different types of trivalent Hopf data
$(B_1,B_2,\mu_l,\nu_l,\mu_r,\nu_r)$ will be taken up with the help of the
table
\raisebox{3pt}{\def\arraystretch{.5}%
\begin{tabular}{|c|c|}\hline\raisebox{-.7pt}{$\sstyle\!\!\!\nu_l\!\!\!$}&
\raisebox{-.7pt}{$\sstyle\!\!\!\nu_r\!\!\!$}\\
\hline\raisebox{-.7pt}{$\sstyle\!\!\!\mu_l\!\!\!$}&
\raisebox{-.7pt}{$\sstyle\!\!\!\mu_r\!\!\!$}\\
\hline\end{tabular}}
where the particular entries take values 0 or 1 dependent on the
(co-)action is trivial or not. Thus by definition maximally three entries
in the table take the value 1.

\begin{corollary}\Label{hp-cases}
\abs\abs
\begin{itemize}
\itembox 0000
All the actions and coactions are trivial. Then the corresponding
Hopf datum is equivalently given by the following data. $B_1$ and $B_2$ are
bialgebras in $\C$, $\Psi_{B_2,B_1}\circ\Psi_{B_1,B_2}=
\id_{B_1\otimes B_2}$, and $B_1\cross{3} B_2$ is the canonical
tensor product bialgebra $B_1\otimes B_2$.
\itembox 1010
The trivalent Hopf datum is given through the following data.
$B_2$ is a bialgebra in $\C$ and $B_1$ is a $B_2$-crossed comodule bialgebra
in ${}^{B_2}_{B_2}\DY{\C}$ \cite{Bes1:97,BD:95}. Then $B_1\cross{3} B_2=
B_1{{}_{\mu_l}^{\nu_l}\!\bowtie} B_2$ is the braided version of the
crossed product or biproduct \cite{Rad1:85,BD:95}.
\itembox 0011
The trivalent Hopf datum $(B_1,B_2,\mu_l,\mu_r)$ is a braided version of
the matched pair \cite{Ma1:90}. Explicitely,
$B_1$ and $B_2$ are bialgebras in $\C$, $B_1$ is a left $B_2$-module
coalgebra, $B_2$ is a right $B_1$-module coalgebra, and the following
defining relations are fulfilled.
\begin{equation}\mathlabel{braid-match-pair}
\begin{array}{c}
\divide\unitlens by 2
\begin{tangle}
\hh\hcu\hstep\id\\
\hh\hstep\ru\\
\hh\hstep\id
\end{tangle}
\multiply\unitlens by 2
=
\divide\unitlens by 3
 \begin{tangle}
      \hh \id\step\cd\step\cd \\
      \id\step\id\step\hx\step\id \\
      \hh \id\step\lu\step\ru \\
      \Ru\dd \\
      \cu
 \end{tangle}
\multiply\unitlens by 3
\quad ,\quad
\mu_r\circ(\eta_2\otimes\id_{B_1})=\eta_2\circ\varepsilon_1\,,\\[20pt]
\divide\unitlens by 2
\begin{tangle}
\hh\id\hstep\hcu\\
\hh\lu\\
\hh\step\id
\end{tangle}
\multiply\unitlens by 2
=
\divide\unitlens by 3
 \begin{tangle}
      \hh \cd\step\cd\step\id \\
      \id\step\hx\step\id\step\id \\
      \hh \lu\step\ru\step\id \\
      \step\d\Lu \\
      \Step\cu
 \end{tangle}
\quad ,\quad
\mu_l\circ(\id_{B_2}\otimes\eta_1)=\eta_1\circ\varepsilon_2\,,\\[20pt]
\divide\unitlens by 3
\begin{tangle}
\hcd\step\hcd\\
\id\step\hx\step\id\\
\ru\step\lu
\end{tangle}
\multiply\unitlens by 3
=
\divide\unitlens by 3
\begin{tangle}
\hcd\step\hcd\\
\id\step\hx\step\id\\
\lu\step\ru\\
\step\hx
\end{tangle}
\multiply\unitlens by 3
\end{array}
\end{equation}
The corresponding bialgebra $B_1 {{}_{\mu_l}\!\bowtie\!{}_{\mu_r}}B_2$
is a braided version of the double cross product \cite{Ma1:90}. The
bialgebra structure reads as follows.
\begin{equation*}
\begin{array}{c}
\m_{B_1 {{}_{\mu_l}\bowtie {}_{\mu_r}} B_2}=
\divide\unitlens by 3
\begin{tangle}
\id\step\hcd\step\hcd\step\id\\
\id\step\id\step\hx\step\id\step\id\\
\id\step\lu\step\ru\step\id\\
\cu\step\cu
\end{tangle}
\multiply\unitlens by 3
\quad ,\quad
\Delta_{B_1 {{}_{\mu_l}\bowtie {}_{\mu_r}} B_2}=
\divide\unitlens by 3
\begin{tangle}
\hcd\step\hcd\\
\id\step\hx\step\id
\end{tangle}
\multiply\unitlens by 3
\end{array}
\end{equation*}
\itembox 1001
The trivalent Hopf datum is given by
$(B_1,B_2,\nu_l,\mu_r)$ where $B_1$ and $B_2$ are bialgebras in $\C$,
$B_1$ is a left $B_2$-comodule coalgebra, $B_2$ is a right $B_1$-module
algebra, and the following defining relations are fulfilled.
\begin{equation*}
\begin{array}{c}
\varepsilon_2\circ\mu_r =\varepsilon_2\otimes\varepsilon_1
\quad ,\quad
\nu_r\circ\eta_2 =\eta_2\otimes\eta_1\,,\\[15pt]
\divide\unitlens by 3
\begin{tangle}
\step\id\step[2]\hcd\\
\step\x\step\id\\
\ld\step[2]\ru\\
\id\step\x\\
\hcu\step[2]\id
\end{tangle}
\multiply\unitlens by 3
=
\divide\unitlens by 3
\begin{tangle}
\id\step\hcd\\
\ru\ld\\
\hcu\step\id
\end{tangle}
\multiply\unitlens by 3
\quad ,\quad
\divide\unitlens by 2
\begin{tangle}
\hh\hstep\ru\\
\hh\hcd
\end{tangle}
\multiply\unitlens by 2
=
\divide\unitlens by 3
 \begin{tangle}
      \hcd\step\cd \\
      \id\step\hx\step\ld \\
      \ru\step\hx\step\id \\
      \cu\step\ru
 \end{tangle}
\multiply\unitlens by 3
\quad ,\quad
\divide\unitlens by 2
\begin{tangle}
\hh\hstep\hcu\\
\hh\ld
\end{tangle}
\multiply\unitlens by 2
=
\divide\unitlens by 3
 \begin{tangle}
      \ld\step\cd \\
      \id\step\hx\step\ld \\
      \ru\step\hx\step\id \\
      \cu\step\hcu
 \end{tangle}\quad .
\multiply\unitlens by 3
\end{array}
\end{equation*}
The affiliated bialgebra $B_1{{}^{\nu_l}\!\bowtie\!{}_{\mu_r}} B_2$
is a braided version of the bicross product \cite{Ma1:90}. The bialgebra
structure is given by
\begin{equation*}
\begin{array}{c}
\m_{B_1{{}^{\nu_l}\bowtie {}_{\mu_r}} B_2}=
\divide\unitlens by 3
\begin{tangle}
\id\step[2]\id\step\hcd\step\id\\
\id\step[2]\hx\step\id\step\id\\
\id\step[2]\id\step\ru\step\id\\
\cu\step\cu
\end{tangle}
\multiply\unitlens by 3
\quad ,\quad
\Delta_{B_1{{}^{\nu_l}\bowtie {}_{\mu_r}} B_2}=
\divide\unitlens by 3
\begin{tangle}
\cd\step\cd\\
\id\step\ld\step\id\step[2]\id\\
\id\step\id\step\hx\step[2]\id\\
\id\step\hcu\step\id\step[2]\id
\end{tangle}
\multiply\unitlens by 3
\end{array}
\end{equation*}
\itembox 1101
This is the most general trivalent Hopf datum $(B_1,B_2,\nu_l,\mu_r,\nu_r)$.
$B_1$ is a bialgebra in $\C$ and a left $B_2$-comodule. $B_2$ is a counital
right $B_1$-module-comodule algebra and a unital coalgebra, and the
following defining identities hold.
\begin{equation*}
\begin{array}{c}
\varepsilon_2\circ\mu_r = \varepsilon_2\otimes\varepsilon_1\,,\quad
\nu_l\circ\eta_1 = \eta_2\otimes\eta_1\\
(\varepsilon_2\otimes\id_{B_1})\circ\nu_r =\eta_1\circ\varepsilon_2\,,
\quad (\id_{B_2}\otimes\varepsilon_1)\circ\nu_l =\eta_2\circ\varepsilon_1
\\[10pt]
\divide\unitlens by 2
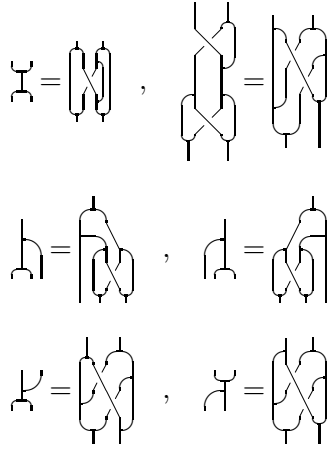
\begin{tangle}
\hh\hcu\\
\hh\hcd
\end{tangle}
\multiply\unitlens by 2
=
\divide\unitlens by 3
\begin{tangle}
\hh\cd\step\cd \\
\hh \id\step\id\step\hrd\hstep\id \\
\id\step\hx\hstep\id\hstep\id \\
\hh \id\step\id\step\hru\hstep\id \\
\hh \cu\step\cu
\end{tangle}
\multiply\unitlens by 3
\quad ,\quad
\divide\unitlens by 3
\begin{tangle}
\step\id\step[2]\hcd\\
\step\x\step\id\\
\step\id\step[2]\ru\\
\ld\step[2]\rd\\
\id\step\x\step\id\\
\hcu\step[2]\hcu
\end{tangle}
\multiply\unitlens by 3
=
\divide\unitlens by 3
\begin{tangle}
\rd\step\cd\\
\id\step\hx\step\ld\\
\id\step\id\step\hx\step\id\\
\ru\step\id\step\hcu\\
\cu\step[1.5]\id
\end{tangle}
\multiply\unitlens by 3
\\ 
\\
\divide\unitlens by 2
\begin{tangle}
\hh\hstep\id\\
\hh\hstep\rd\\
\hh\hcd\hstep\id
\end{tangle}
\multiply\unitlens by 2
=
\divide\unitlens by 3
\begin{tangle}
\cd \\
\Rd\d \\
\hh \id\step\ld\step\rd \\
\id\step\id\step\hx\step\id \\
\hh \id\step\cu\step\cu
\end{tangle}
\multiply\unitlens by 3
\quad ,\quad
\divide\unitlens by 2
\begin{tangle}
\hh\step\id\\
\hh\ld\\
\hh\id\hstep\hcd
\end{tangle}
\multiply\unitlens by 2
=
\divide\unitlens by 3
\begin{tangle}
\Step\cd \\
\step\dd\Ld \\
\hh \ld\step\rd\step\id \\
\id\step\hx\step\id\step\id \\
\hh \cu\step\cu\step\id
\end{tangle}
\multiply\unitlens by 3 
\\ 
\\
\divide\unitlens by 2
\begin{tangle}
\hh\hstep\ru\\
\hh\hcd
\end{tangle}
\multiply\unitlens by 2
=
\divide\unitlens by 3
\begin{tangle}
\hcd\step\cd \\
\id\step\hx\step\ld \\
\ru\step\hx\step\id \\
\cu\step\ru
\end{tangle}
\multiply\unitlens by 3
\quad ,\quad
\divide\unitlens by 2
\begin{tangle}
\hh\hstep\hcu\\
\hh\ld
\end{tangle}
\multiply\unitlens by 2
=
\divide\unitlens by 3
\begin{tangle}
\ld\step\cd \\
\id\step\hx\step\ld \\
\ru\step\hx\step\id \\
\cu\step\hcu
\end{tangle}
\multiply\unitlens by 3
\end{array}
\end{equation*}
The structure of the resulting bialgebra
$B_1{{}^{\nu_l}\!\bowtie\!{}^{\nu_r}_{\mu_r}} B_2$ is given by
\begin{equation*}
\begin{array}{c}
\m_{B_1{{}^{\nu_l}\bowtie {}^{\nu_r}_{\mu_r}} B_2}=
\divide\unitlens by 3
\begin{tangle}
\id\step[2]\id\step\hcd\step\id\\
\id\step[2]\hx\step\id\step\id\\
\id\step[2]\id\step\ru\step\id\\
\cu\step\cu
\end{tangle}
\multiply\unitlens by 3
\quad ,\quad
\Delta_{B_1{{}^{\nu_l}\bowtie {}^{\nu_r}_{\mu_r}} B_2}=
\divide\unitlens by 3
\begin{tangle}
\cd\step\cd\\
\id\step\ld\step\rd\step\id\\
\id\step\id\step\hx\step\id\step\id\\
\id\step\hcu\step\hcu\step\id
\end{tangle}
\multiply\unitlens by 3
\end{array}
\end{equation*}
\end{itemize}
\end{corollary}

\begin{proof}
Straightforward evaluations of Definition \ref{hopf-pair} using the
particular trivial (co-)\-actions.
\end{proof}
\abs
In the following we will discuss two examples which are closely
related to Hopf algebras constructed by Ore extensions \cite{BDG:98}.
The first example is an infinite Hopf algebra whereas the second
example is Radford's finite-dimensional 4-parameter Hopf algebra
\cite{Rad1:94}. Both Hopf algebra types turn out to be trivalent cross
product bialgebras of type \crossbox 0101 and therefore are biproduct
Hopf algebras according to \cite{Rad1:85}. Closely related are the examples 
studied in the forthcoming paper \cite{CIMZ:98} where the universal
aspects have been considered.

\begin{example}\Label{ex-semi-ore}
{\normalfont Suppose $C$ is an abelian group, $k$ is an algebraically
closed field with $\mathrm{char}(k)=0$. Let $C^*$ be the chartacter
group of $C$ and let $t\in\NN$. Assume further that 
$g=(g_1,\ldots,g_t)\in C^t$ and $g^*=(g^*_1,\ldots,g^*_t)\in (C^*)^t$
where at least one of the $g_i$ and one of the $g^*_j$ is non-trivial,
and $g_{lr}:= g^*_l(g_r)$ with $g_{lr}\cdot g_{rl}=1$ for any 
$l,r\in\{1,\ldots, t\}$. The algebra
$H(C,t,g,g^*)$ be generated by the group $C$ and the
generators $\{x_i\}_{i=1}^t$ subject to the additional relations
\begin{equation}\mathlabel{gen-rels1}
x_j\cdot c = g_j^*(c)\,c\cdot x_j\quad\text{and}\quad
x_j\cdot x_k = g_{jk}\,x_k\cdot x_j
\end{equation}
for any $j,k\in\{1,\ldots,t\}$ and $c\in C$.
Then the following relations define a Hopf algebra structure on
$H(C,t,g,g^*)$.
\begin{equation}\mathlabel{gen-hopf1}
{\begin{gathered}
\Delta(c)=c\otimes c\,,\\
\Delta(x_j)=x_j\otimes g_j +\E\otimes x_j\,,
\end{gathered}}
\quad
{\begin{gathered}
\varepsilon(c)=1\,,\\
\varepsilon(x_j)=0\,,
\end{gathered}}
\quad
{\begin{gathered}
S(c)=c^{-1}\,,\\
S(x_j)=-x_j\cdot g_j^{-1}
\end{gathered}}
\end{equation}
for all $j,k\in\{1,\ldots,t\}$ and $c\in C$.
Every element of $H(C,t,g,g^*)$ is a finite sum of the form
$h=\sum_{c\in C} c\cdot f(x_j)_c$ where $f(x_j)_c$ is a 
(non-commutative) polynomial in $\{x_j\}$. The Hopf algebra is
finite-dimensional if $C$ is finite and $g_{jj}=-1$ for all 
$j\in\{1,\ldots,t\}$. 

Now we define $B_1:= kC$ the group Hopf algebra of $C$ with
comultiplication $\Delta_1$ and counit $\varepsilon_1$ on $C$. Let 
$B_2$ be the algebra
$B_2:= k\langle\{x_j\}\rangle / (x_j\cdot x_k = g_{jk}\,x_k\cdot x_j
\ \forall j,k\in\{1,\ldots,t\})$. Then the following definitions
yield $k$-linear morphisms.
\begin{equation}\mathlabel{inj-proj1}
{\begin{split}
\inj_1 &:\left\{\begin{matrix} B_1\hookrightarrow H(C,t,g,g^*)\\
c\mapsto c\end{matrix}\right.\\
\proj_1&:\left\{\begin{matrix}H(C,t,g,g^*)\to B_1\\
c\cdot f(x_j)\mapsto \varepsilon(f(x_j))\,c\end{matrix}\right.
\end{split}}
\quad
{\begin{split}
\inj_2&:\left\{\begin{matrix}B_2\hookrightarrow H(C,t,g,g^*)\\
x\mapsto x\end{matrix}\right.\\
\proj_2&:\left\{\begin{matrix}H(C,t,g,g^*)\to B_2\\
c\cdot f(x_j)\mapsto \varepsilon(c)\,f(x_j)\end{matrix}\right.
\end{split}}
\end{equation}
Straightforward calculations show that $\inj_1$ is a Hopf algebra
morphism, $\inj_2$ is an algebra morphism, and $\proj_1$ is an algebra
morphism. Then one concludes easily that $\proj_1$ is coalgebra
morphism since $(\Delta_1\otimes\Delta_1)\circ\proj_1$ and
$\proj_1\circ\Delta$ are algebra morphism, and the equality of both
morphisms has to be proven on the generators only. 
According to \cite{CIMZ:98} we define
the comultiplication $\Delta_2$ and counit $\varepsilon_2$ of $B_2$ 
as $\Delta_2:= (\proj_2\otimes\proj_2)\circ\Delta\circ\inj_2$ and
$\varepsilon_2:= \varepsilon\circ\inj_2$. 
Using again the fact that every element of $H(C,t,g,g^*)$ 
is a finite sum of the form $h=\sum_{c\in C} c\cdot f(x_j)_c$, one finds
\begin{equation*}
\begin{split}
\Delta_2\circ\proj_2(h)&=
\sum_{c\in C} \varepsilon(c)\, \Delta_2(f(x_j)_c)\\
&=\sum_{c\in C} \Delta(f(x_j)_c)\\
&=\sum_{c\in C} (\proj_2\otimes\proj_2)\circ\Delta(f(x_j)_c)
\end{split}
\end{equation*}
and on the other hand
\begin{equation*}
\begin{split}
(\proj_2\otimes\proj_2)\circ\Delta(h)&=\sum_{c\in C} (\proj_2\otimes\proj_2)
\big((c\otimes c)\cdot \Delta(f(x_j)_c)\big)\\
&=\sum_{c\in C} (\proj_2\otimes\proj_2)\big(\Delta(f(x_j)_c)\big)\,.
\end{split}
\end{equation*}
Hence $\Delta_2\circ\proj_2= (\proj_2\otimes\proj_2)\circ\Delta$. 
Then it follows that $\Delta_2$ is coassociative and $\proj_2$ is 
coalgebra morphism. Furthermore $\proj_2$ is not an algebra morphism.
Suppose the converse, then 
$\proj_2(c\cdot x_j)=\proj_2(c)\cdot\proj_2(x_j)=
\varepsilon(c)\,x_j$. On the other hand $\proj_2(c\cdot x_j)=
\proj_2(g_j^*(c)\,x_j\cdot c)= g_j^*(c)\varepsilon(c)\,x_j$. 
Therefore $\varepsilon(c)\,x_j= g_j^*(c)\varepsilon(c)\,x_j$. But by
assumption there exists a non-trivial $g_j^*$ which then leads to a
contradiction. Similarly, by the existence of a non-trivial $g_j$ it
can be shown that $\inj_2$ is not a coalgebra morphism.

Now we prove that $(\proj_1\otimes\proj_2)\circ\Delta =
\big(\m\circ(\inj_1\otimes\inj_2)\big)^{-1}$.
Observe that $(\proj_1\otimes \id)\circ\Delta\circ\inj_2=
(\E\otimes \id)\circ\inj_2$ 
because the corresponding identity holds on the generators 
and then the statement follows since 
$(\proj_1\otimes \id)\circ\Delta\circ\inj_2$ and $(\E\otimes\id)\circ\inj_2$
are algebra morphisms. Then for any finite sum 
$\sum_{c\in C} c\otimes f(x_j)_c \in B_1\otimes B_2$ we have
\begin{equation*}
\begin{split}
(\proj_1\otimes\proj_2)\circ\Delta\circ\m\circ
(\inj_1\otimes\inj_2)\sum_{c\in C} c\otimes f(x_j)_c &=
\sum_{c\in C}(\proj_1\otimes\proj_2)\circ\Delta(c\cdot f(x_j)_c)\\
&=\sum_{c\in C}(c\otimes\E)\cdot(\proj_1\otimes\proj_2)\circ
\Delta(f(x_j)_c)\\
&=\sum_{c\in C}(c\otimes\E)\cdot(\E\otimes\proj_2(f(x_j)_c))\\
&=\sum_{c\in C}c\otimes f(x_j)_c\,.
\end{split}
\end{equation*}
Thus $(\proj_1\otimes\proj_2)\circ\Delta$ is a left inverse of
$\m\circ(\inj_1\otimes\inj_2)$. Similarly the right invertibility
can be proven. Therefore all conditions of Theorem \ref{rhp-cpb}
hold, and the morphisms $\varphi_{1,2}$ and $\varphi_{2,1}$ of the
corresponding cross product bialgebra are given by
$\varphi_{1,2} = (\varepsilon_1\otimes\id_{B_2}\otimes
\id_{B_1}\otimes\varepsilon_2)\circ\Delta_B$ and
$\varphi_{2,1}=\m_B\circ(\eta_1\otimes\id_{B_2}\otimes\id_{B_1}\otimes
\eta_2)$ (see the proof of Proposition \ref{crossprod-proj}).
Using \eqref{bat-act} yields the (non-trivial) (co-)actions 
\begin{equation*}
\mu_r(x\otimes c) = \proj_2(x\cdot c)\quad\text{and}\quad
\nu_r(x) = (\proj_2\otimes\proj_1)\circ\Delta(x)\,.
\end{equation*}
From Corollary \ref{hp-cases} if follows eventually
\begin{proposition}
The Hopf algebra $H(C,t,g,g^*)$ is isomorphic to
the biproduct bialgebra $B_1\bowtie\!_{\mu_r}^{\nu_r} B_2$.
In particular $B_2$ is a right $B_1$-crossed module bialgebra.\endproof
\end{proposition}}
\end{example}
\abs
The following example is Radford's 4-parameter Hopf algebra \cite{Rad1:94}.
Its structure resembles the one of the Hopf algebra 
$H(C,t,g,g^*)$ discussed in Example \ref{ex-semi-ore}. 

\begin{example}[Radford's 4-Parameter Hopf Algebra]\Label{ex-rad}
{\normalfont Suppose again that $k$ is an algebraically closed field
with $\mathrm{char}(k)=0$. Let $n,N,\nu$ be positive integers such
that $n\vert N$, and $\nu < n$. Let $q$ be a primitive $n$.th root
of unity and $q^\nu$ be an $r$.th root of unity, where $r=n/(n,\nu)$.
Set $C_N$ to be the cyclic group of order $N$ and $g$ a generating
element of $C_N$. Then the 4-parameter Hopf algebra $H_{n,q,N,\nu}$
is generated by the group algebra $k C_N$ and the generator $x$ subject
to the additional relations
\begin{equation}\mathlabel{gen-rels2}
x^r =0\quad\text{and}\quad x\cdot g=q\,g\cdot x\,.
\end{equation}
On the generators the comultiplication $\Delta$, counit $\varepsilon$ and 
antipode $S$ are given by
\begin{equation*}
{\begin{gathered}
\Delta(g)=g\otimes g\,,\\
\Delta(x)=x\otimes\E+g^{-\nu}\otimes x\,,
\end{gathered}}
\quad
{\begin{gathered}
\varepsilon(g)=1\,,\\
\varepsilon(x)=0\,,
\end{gathered}}
\quad
{\begin{gathered}
S(g)=g^{-1}\,,\\
S(x)=-g^\nu\cdot x\,,
\end{gathered}}
\end{equation*}
Similarly as in the previous example one shows that every element of
$H_{n,q,N,\nu}$ can be represented in the form $b=\sum_{m=0,l=0}^{r-1,N-1}
\lambda_{m,l}\,x^m\cdot g^l$. In particular $H_{n,q,N,\nu}$ is
finite-dimensional. Let $B_1$ be the algebra $B_1 := k\langle x\rangle/(x^r)$ 
and $B_2$ be the group Hopf algebra $B_2:=k C_N$.
Then $B_1$ becomes a coalgebra through
\begin{equation}\mathlabel{gen-hopf2}
\Delta_1(x^m)=
\sum_{l=0}^m \left(\begin{smallmatrix}m\\l\end{smallmatrix}\right)_{q^\nu}
x^l\otimes x^{m-l} \quad\text{and}\quad \varepsilon_1(x^m)=\delta_{m,0}
\end{equation}
where the $q$-binomial is 
$\left(\begin{smallmatrix}m\\l\end{smallmatrix}\right)_{p}:=
\frac{(m)_p!}{(l)_p!\,(m-l)_p!}$, $(s)_p!:=
(1)_p\cdot(2)_p\cdot\ldots\cdot(s)_p$, $(0)_p :=1$, and 
$(s)_p := \frac{1-p^s}{1-p}$. Using the properties of the
$q$-binomials it is evident that \eqref{gen-hopf2} renders
$B_1$ a coalgebra. Then it can be proven that
the following definitions yield $k$-linear homomorphisms 
$\inj_1$, $\inj_2$, $\proj_1$ and $\proj_2$.
\begin{equation}\mathlabel{inj-proj2}
{\begin{split}
\inj_1&:\left\{\begin{matrix}B_1\hookrightarrow H_{n,q,N,\nu}\\
x^m\mapsto x^m\end{matrix}\right.\\
\proj_1&:\left\{\begin{matrix}H_{n,q,N,\nu}\to B_1\\
x^m\cdot g^l\mapsto x^m\end{matrix}\right.
\end{split}}
\quad
{\begin{split}
\inj_2&:\left\{\begin{matrix}B_2\hookrightarrow H_{n,q,N,\nu}\\
g^l\mapsto g^l\end{matrix}\right.\\
\proj_2&:\left\{\begin{matrix}H(C,t,g,g^*)\to B_2\\
x^m\cdot g^l\mapsto \delta_{m,0}\,g^l\end{matrix}\right.
\end{split}}
\end{equation}
Obviously $\inj_1$, $\inj_2$ and $\proj_2$ are
algebra morphisms since they preserve the relations of the
algebras $H_{n,q,N,\nu}$, $B_1$ and $B_2$. Furthermore
$\inj_2$ and $\proj_2$ are coalgebra homomorphisms since
the corresponding identities hold for the generators $g$
and $x$. Finally, $\proj_1$ is coalgebra morphism because
$\Delta_1\circ\proj_1(x^m\cdot g^l)=\Delta_1(x^m)$ and
\begin{equation*}
\begin{split}
(\proj_1\otimes\proj_1)\circ\Delta(x^m\cdot g^l)&=
(\proj_1\otimes\proj_1)\big(\Delta(x)^m\cdot (g^l\otimes g^l)\big)\\
&=(\proj_1\otimes\proj_1)\big((x\otimes \E + g^{-\nu}\otimes x)^m\cdot 
(g^l\otimes g^l)\big)\\
&=\sum_{j=0}^m \left(\begin{smallmatrix}m\\j\end{smallmatrix}\right)_{q^\nu}\,
(\proj_1\otimes\proj_1)(x^j\cdot g^{l-\nu\,(m-j)}\otimes x^{m-j}\cdot g^l)\\
&=\sum_{j=0}^m\left(\begin{smallmatrix}m\\j\end{smallmatrix}\right)_{q^\nu}
\,x^j\otimes x^{m-j}\\
&=\Delta_1(x^m)
\end{split}
\end{equation*}
where we used the $q$-binomial identity $(a+b)^m =\sum_{j=0}^m
\left(\begin{smallmatrix}m\\j\end{smallmatrix}\right)_\lambda a^j\cdot b^{m-j}$
if $a\cdot b= \lambda^{-1}\,b\cdot a$. This proves that $\proj_1$ is
a coalgebra morphism. Since by assumption $q\ne 1$ and $g^\nu\ne\E$
one concludes similarly as in Example \ref{ex-semi-ore} that
$\inj_1$ in no coalgebra morphism and $\proj_1$ is no algebra
morphism. Because $H_{n,q,N,\nu}$ is finite-dimensional the subsequent
relations prove that $(\proj_1\otimes\proj_2)\circ\Delta = 
\big(\m\circ(\inj_1\otimes\inj_2)\big)^{-1}$.
\begin{equation*}
\begin{split}
(\proj_1\otimes\proj_2)\circ\Delta\circ\m\circ
(\inj_1\otimes\inj_2)(x^m\otimes g^l) &=
(\proj_1\otimes\proj_2)\circ\Delta (x^m\cdot g^l)\\
&=(\proj_1\otimes\proj_2)\big(\Delta(x)^m\cdot (g^l\otimes g^l)\big)\\
&=\sum_{j=0}^m \left(\begin{smallmatrix}m\\j\end{smallmatrix}\right)_{q^\nu}\,
(\proj_1\otimes\proj_2)(x^j\cdot g^{l-\nu\,(m-j)}\otimes x^{m-j}\cdot g^l)\\
&=\sum_{j=0}^m \left(\begin{smallmatrix}m\\j\end{smallmatrix}\right)_{q^\nu}\,
x^j\otimes \varepsilon(x^{m-j})\,g^l)\\
&= x^m\otimes g^l\,.
\end{split}
\end{equation*}
Hence again all conditions of Theorem \ref{rhp-cpb}
hold, and the morphisms $\varphi_{1,2}$ and $\varphi_{2,1}$ of the
corresponding cross product bialgebra are given by
$\varphi_{1,2} = (\varepsilon_1\otimes\id_{B_2}\otimes
\id_{B_1}\otimes\varepsilon_2)\circ\Delta_B$ and
$\varphi_{2,1}=\m_B\circ(\eta_1\otimes\id_{B_2}\otimes\id_{B_1}\otimes
\eta_2)$, or explicitely $\varphi_{1,2}(x^m\otimes g^l)=g^{l-\nu m}\otimes
x^m$ and $\varphi_{2,1}(g^l\otimes x^m)= q^{-ml}\,x^m\otimes g^l$.
Using \eqref{bat-act} yields the (non-trivial) (co-)actions 
\begin{equation*}
\mu_l(g^l\otimes x^m) = q^{-ml}\,x^m \quad\text{and}\quad
\nu_l(x^m)= g^{-\nu m}\otimes x^m\,.
\end{equation*}
We collect the previous results in the next proposition.
\begin{proposition}
The 4-parameter Hopf algebra $H_{n,q,N,\nu}$ is isomorphic to
the biproduct bialgebra $B_1{}_{\mu_r}^{\nu_r}\!\bowtie B_2$.
In particular $B_1$ is a left $B_2$-crossed module bialgebra.\endproof
\end{proposition}}
\end{example}

\section{Applications}\Label{applic}

In Section \ref{applic} we discuss further applications of the results
of Section \ref{cross-prod}. The first example shows that Majid's
double biproduct \cite{Ma1:95} essentially is a cross product 
bialgebra construction in the braided category $\hhchh$
of Hopf bimodules over a Hopf algebra $H$.
In the second subsection we show that the Drinfel'd double in a 
braided category can be reconstructed as a matched pair if and only 
if the braiding of the two tensor factors is involutive. This confirms
in a certain sense the statements of \cite{Ma1:95}.

\subsection{Double Biproduct Bialgebras}

We consider a Hopf algebra $H$. From Theorem \ref{yd-hopfbi} and its mirror
symmetric version we know that every
right $H$-crossed module $B\in \Obj\big(\DY\C^H_H\big)$ and every left
$H$-crossed module $C\in\Obj\left({}^H_H\DY\C\right)$ yield
Hopf bimodules $X= H\ltimes B$ and $Y=C\rtimes H$ in $\hhchh$
respectively. These special types of Hopf bimodules
will be used later for the construction of the double biproduct
as a twist of a tensor product bialgebra in $\hhchh$ considered as
bialgebra in the category $\C$ according to Proposition \ref{hbb-bp}.

\begin{lemma}\Label{psi-2}
Let $B\in\Obj\big(\DY\C^H_H\big)$ and
$C\in\Obj\big({}^H_H\DY\C\big)$ be $H$-crossed modules.
Then for the objects $X=H\ltimes B$ and $Y=C\rtimes H$
in the category $\hhchh$ the identity
${}^\hhchh\Psi_{X,Y}\circ{}^\hhchh\Psi_{Y,X}=\id_{X\otimes_HY}$ holds if
and only if
\begin{equation}
\mathlabel{hhchhpsi2triv}
\Psi_{C,B}\circ\Psi_{B,C}=(\mu^B_r\otimes\mu^C_l)\circ
 (\id_B\otimes\Psi_{H,H}\otimes\id_C)\circ (\nu^B_r\otimes\nu^C_l)\,.
\end{equation}
\end{lemma}

\begin{proof}
From \eqref{wor-braid1} and its dual version one deduces
\begin{equation}\mathlabel{hhchhpsi2}
\begin{split}
&({}_X\Pi\otimes\Pi_Y)\circ\rho_{X,Y}\circ{}^\hhchh\Psi_{Y,X}\circ
 {}^\hhchh\Psi_{X,Y}\circ\lambda_{X,Y}\circ({}_X\Pi\otimes\Pi_Y)\\
&=({}_X\Pi\otimes \Pi_Y)\circ\Psi_{Y,X}\circ\rho^H_{Y,X}\circ
 \lambda^H_{Y,X}\circ\Psi_{X,Y}\circ({}_X\Pi\otimes\Pi_Y)\,.
\end{split}
\end{equation}
Condition \eqref{hhchhpsi2triv} means that the right hand side
of \eqref{hhchhpsi2} equals to
$({}_X\Pi\otimes\Pi_Y)\circ\rho_{X,Y}\circ\lambda_{X,Y}\circ
({}_X\Pi\otimes\Pi_Y)$. Then we use Lemma \ref{aux-psi-inv} to derive
$\rho_{X,Y}\circ{}^\hhchh\Psi_{X,Y}^2\circ\lambda_{X,Y}
 =\rho_{X,Y}\circ\lambda_{X,Y}$.
Since $\lambda^H_{X,Y}$ is an epimorphism and $\rho^H_{X,Y}$ is a
monomorphism the sufficient part of the lemma is proved.
Conversely if ${}^\hhchh\Psi_{X,Y}\circ{}^\hhchh\Psi_{Y,X}=
\id_{X\otimes_HY}$ then equation \eqref{rho-lambda} proves
\eqref{hhchhpsi2triv}.\end{proof}
\abs
Lemma \ref{psi-2} is a braided version of the identity (58)
in \cite{Ma1:95} which was one of the compatibility conditions
for the construction of the double biproduct bialgebra.
We suppose henceforth that $B$ and $C$ are bialgebras in $\DY\C^H_H$ and
${}^H_H\DY\C$ respectively. Then according to Theorem \ref{yd-hopfbi} the
objects $X=H\ltimes B$ and $Y=C\rtimes H$ are bialgebras in $\hhchh$.
From \cite{BD:95} (and Proposition \ref{hbb-bp}) we know that the
multiplications $\underline\m_X$ and $\underline\m_Y$ are uniquely
determined. In particular $\m_{H\ltimes B}=
\underline\m_X\circ\lambda^H_{X,X}$ and $\underline\m_X=\m_{H\ltimes B}
\circ(\id_H\otimes{}_X\inj)$ where $\m_{H\ltimes B}$ is the
multiplication of the crossed product algebra in $\C$
(see Corollary \ref{hp-cases} and \cite{BD:95}) and ${}_X\inj$ is the
morphism described below Theorem \ref{Hopf-br}. Similarly the
multiplication $\underline\m_Y$ and the comultiplications
$\underline\Delta_X$ and $\underline\Delta_Y$ can be calculated. The
condition ${}^\hhchh\Psi_{X,Y}\circ{}^\hhchh\Psi_{Y,X}=\id_{Y\otimes_H X}$
of Lemma \ref{psi-2} allows the construction of the canonical tensor
product bialgebra $Y\otimes_HX$ in $\hhchh$ according to Corollary
\ref{hp-cases}. Explicitely
\begin{equation}\mathlabel{mult-comult-hobi}
\begin{split}
\m_{Y\otimes_HX}&=(\underline\m_Y\otimes_H\underline\m_X)
\circ(\id_Y\otimes_H{}^\hhchh\Psi_{X,Y}\otimes_H\id_X)\,,\\
\Delta_{Y\otimes_HX}&=
(\id_Y\otimes_H{}^\hhchh\Psi_{Y,X}\otimes_H\id_X)\circ
(\underline\Delta_Y\otimes_H\underline\Delta_X)\,.
\end{split}
\end{equation}

\begin{lemma}\Label{chb-bialg}
Suppose that the conditions of Lemma \ref{psi-2} are fulfilled for the
bialgebras $X$ and $Y$ in $\hhchh$. Then the Hopf bimodule tensor product
bialgebra $Y\otimes_HX$ is identified with the object
$C\otimes H\otimes B$ through the canonical morphisms given by
\begin{equation}
\mathlabel{lr-special}
\lambda^H_{Y,X}=\id_C\otimes\m_H\otimes\id_B,\qquad
\rho^H_{Y,X}=\id_C\otimes\Delta_H\otimes\id_B\,.
\end{equation}
We denote by $\underline Z := Y\otimes_H X=C\otimes H\otimes B$ the bialgebra 
in $\hhchh$. Then
\begin{equation}\mathlabel{lr-special1}
\lambda^H_{\underline Z,\underline Z}=\id_C\otimes\lambda^H_{X,Y}\otimes
 \id_B\,,\quad\rho^H_{\underline Z,\underline Z} =
 \id_C\otimes\rho^H_{X,Y}\otimes\id_B\,.
\end{equation}
\end{lemma}

\begin{proof}
Without problems one verifies that $\lambda^H_{Y,X}$ and
$\rho^H_{Y,X}$ according to \eqref{lr-special} fulfill the required
properties of Theorem \ref{Hopf-br}. Using the identification
$Y\otimes_H X\cong C\otimes H\otimes B$ induced from \eqref{lr-special}
the proof of \eqref{lr-special1} follows.
\end{proof}
\abs
According to Proposition \ref{hbb-bp} one obtains a bialgebra
$Z=C\otimes H\otimes B$ in $\C$ from the $H$-Hopf bimodule bialgebra
$\underline Z = Y\otimes_H X$. The explicit structure
of $Z$ is presented in the subsequent proposition.

\begin{proposition}\Label{free-biprod}
Multiplication and comultiplication of the bialgebra
$Z=C\otimes H\otimes B$ are given by
\begin{equation}\mathlabel{mz}
\m_Z=
\begin{array}{c}
\divide\unitlens by 2
\begin{tangle}
 \object{\sstyle C}\step\hstep\object{\sstyle H}\step\hstep%
 \object{\sstyle B}\step\object{\sstyle C}\step\hstep\object{\sstyle H}%
 \step\hstep\object{\sstyle B}\\
 \id\step\hcd\step\hx\step\hcd\step\id\\
 \id\step\id\step\hx\step\hx\step\id\step\id \\
 \hh \id\step\lu\step\cu\step\ru\step\id \\
 \cu\step\hstep\id\hstep\step\cu\\
\step\object{\sstyle C}\step[2.5]\object{\sstyle H}\step[2.5]%
 \object{\sstyle B}\step
\end{tangle}
\multiply\unitlens by 2
\end{array}
\quad\text{and}\quad
\Delta_Z =
 \begin{array}{c}
\divide\unitlens by 2
\begin{tangle}
 \step\object{\sstyle C}\step[2.5]\object{\sstyle H}\step[2.5]%
 \object{\sstyle B}\step\\
 \cd\step\hstep\id\hstep\step\cd \\
 \hh \id\step\ld\step\cd\step\rd\step\id \\
 \id\step\id\step\hx\step\hx\step\id\step\id \\
 \id\step\hcu\step\hx\step\hcu\step\id\\
 \object{\sstyle C}\step\hstep\object{\sstyle H}\step\hstep%
 \object{\sstyle B}\step\object{\sstyle C}\step\hstep\object{\sstyle H}%
 \step\hstep\object{\sstyle B}
\end{tangle}
\multiply\unitlens by 2
\end{array}
\end{equation}
The unit and counit read as $\eta_Z=\eta_C\otimes\eta_H\otimes\eta_B$ and
$\varepsilon_Z=\varepsilon_C\otimes\varepsilon_H
\otimes\varepsilon_B$ respectively.
There are canonical bialgebra monomorphisms
$\id_{C\otimes H}\otimes\eta_B:C\rtimes H\to Z$,
$\eta_C\otimes\id_{H\otimes B}:H\ltimes B\to Z$, and
bialgebra epimorphisms
$\id_{C\otimes H}\otimes\varepsilon_B:Z\to
C\rtimes H$, $\varepsilon_C\otimes\id_{H\otimes B}:Z
\to H\ltimes B$. Additionally
$\id_C\otimes\eta_H\otimes\id_B:C\otimes B\to Z$ is an algebra
monomorphism, and $\id_C\otimes\varepsilon_H\otimes\id_B:Z\to C\otimes B$
is a coalgebra epimorphism.
\end{proposition}

\begin{proof}
From the considerations before Lemma \ref{chb-bialg} we know that
the multiplication of $\underline Z=Y\otimes_H X$ is given by
$\m_{\underline Z}=(\underline\m_Y\otimes_H\underline\m_X)\circ
 (\id_Y\otimes_H{}^\scripthhchh\Psi_{X,Y}\otimes_H\id_X)$.
and according to Proposition \ref{hbb-bp} the multiplication of $Z$ reads
as
\begin{equation}\mathlabel{mz1}
\begin{split}
\m_Z &=\m_{\underline Z}\circ\lambda^H_{\underline Z,\underline Z}\\
&=(\m_Y\otimes_H\m_X)\circ
  (\id_Y\otimes_H\otimes{}^\hhchh\Psi_{X,Y}\otimes_H\id_X)\circ
  \lambda^H_{\underline Z,\underline Z}\\
&=(\m_C\otimes H\otimes\m_B)\circ
  (\id_C\otimes{}^\hhchh\Psi_{X,Y}\circ\lambda^H_{X,Y}\otimes\id_B)
\end{split}
\end{equation}
where we used \eqref{lr-special1} in the third equation.
From \eqref{wor-braid1} and from the Hopf bimodule property of
${}^\hhchh\Psi_{X,Y}$ and $\lambda^H_{X,Y}$ (see Theorem \ref{Hopf-br})
it follows
\begin{equation}\mathlabel{mz2}
\begin{split}
&{}^\hhchh\Psi_{X,Y}\circ\lambda^H_{X,Y}\\
&={}^\hhchh\Psi_{X,Y}\circ\lambda^H_{B,C}\circ
 \left(\left(\mu^X_l\circ(\id_H\otimes{}_X\Pi)\circ\nu^X_l\right)\otimes
 \left(\mu^Y_r\circ(\Pi_Y\otimes\id_H)\circ\nu^Y_r\right)\right)\\
&=M\circ\left(\id_H\otimes\left(\lambda^H_{Y,X}\circ\Psi_{X,Y}
 \circ({}_X\Pi\otimes\Pi_Y)\right)\otimes\id_H\right)
 \circ(\nu^X_l\otimes\nu^Y_r)\\
&=\lambda^H_{Y,X}\circ(\mu^Y_\ell\otimes\mu^X_r)\circ
   (\id_H\otimes\Psi_{X,Y}\otimes\id_H)
\end{split}   
\end{equation}
where $M=\mu^{Y\otimes_HX}_r\circ(\mu^{Y\otimes_HX}_l\otimes\id_H)$.
Inserting \eqref{mz2} into \eqref{mz1} yields the final result for
$\m_Z$. Dually $\Delta_Z$ can be derived.
\end{proof}
\abs
If it is clear from the context we will henceforth denote the bialgebra
$Z$ in $\C$ by $C\otimes H\otimes B$.

\begin{proposition}\Label{dbiprod}
Let $B$ and $C$ be as in Lemma \ref{chb-bialg}.
Suppose that $\rho:B\otimes C\to\1_\C$ is a morphism in $\C$ satisfying
the identities
\begin{equation}\mathlabel{dbiprod-coc}
\begin{split}
\rho\circ(\mu^B_r\otimes\id_B)&=\rho\circ(\id_B\otimes\mu^C_\ell)\,,\\
\rho\circ(\id_B\otimes\m_C)&=
 \rho^{(2)}\circ(\Psi^{-1}_{B,B}\circ\Delta_B\otimes\id_{C\otimes C})\,,\\
\rho\circ(\m_B\otimes\id_C)&=\rho^{(2)}\circ({}^{\DY\C^H_H}\Psi_{B,B}
 \otimes\Psi^{-1}_{C,C}\circ\Delta_C)
\end{split}
\end{equation}
where $\rho^{(2)}=\rho\circ(\id_B\otimes\rho\otimes\id_C)$ and
${}^{\DY\C^H_H}\Psi_{B,B}=(\id_B\otimes\mu_r^B)\circ(\Psi_{B,B}\otimes
\id_H)\circ(\id_B\otimes\nu_r^B)$ is the braiding in $\DY\C^H_H$
\cite{BD:95}. Then $\hat\rho:=\varepsilon_{C\rtimes H}\otimes\rho\otimes
\varepsilon_{H\ltimes B}$ is a
$2$-cocycle of the bialgebra $C\otimes H\otimes B$ from Proposition
\ref{free-biprod}. If $\rho$ is
convolution invertible then $\hat\rho$ is convolution invertible with
inverse $\hat\rho^{\text{--}}= \varepsilon_{C\rtimes H}\otimes
\rho^{\text{--}}\otimes\varepsilon_{H\ltimes B}$.
\end{proposition}

\begin{proof}
The convolution invertibility and the cocycle property \eqref{2cocycle2}
for $\hat\rho$, $\m_Z$ and $\Delta_Z$ can be proven easily. The
calculation of the left and the right hand side of \eqref{2cocycle1}
respectively yields
\begin{equation}\mathlabel{2cocycle3}
\begin{split}
&\rho\circ(\mu_r^B\otimes\id_C)\circ\big((\id_B\otimes\rho)\circ(\Psi_{B,B}
\otimes\id_C)\circ(\Delta_B\otimes\id_C)\otimes\\
&\otimes((\rho\otimes\id_H)\circ(\id_B\otimes\Psi_{H,C})\otimes\id_C)\circ\\
&\circ((\id_B\otimes\m_H)\circ(\Psi_{H,B}\otimes\id_H)\circ
(\id_H\otimes\nu_r^B)
\otimes\Delta_C)\big)
\end{split}
\end{equation}
and
\begin{equation}\mathlabel{2cocycle4}
\begin{split}
&\rho\circ(\m_B\otimes\id_C)\circ(\mu_r^B\otimes\id_B\otimes\id_C)\circ\\
&\circ\big((\id_B\otimes\rho)\circ(\Pi_{B,B}\circ\Delta_B\otimes\id_C)
\otimes\id_{H\otimes B\otimes C}\big)
\end{split}
\end{equation}
where the first and the second defining property of $\rho$ in
\eqref{dbiprod-coc} have been used. With the help of the third
equation of \eqref{dbiprod-coc} one can show that \eqref{2cocycle3}
equals \eqref{2cocycle4}.\end{proof}
\abs
\begin{remark}
{\normalfont The identities \eqref{dbiprod-coc} are braided versions of
\cite[eqs.~(56)]{Ma1:95}.}
\end{remark}
\abs
The following corollary is a straightforward consequence of
Propositions \ref{twist-bialg}, \ref{free-biprod} and \ref{dbiprod}.

\begin{corollary}\Label{double-biprod}
Suppose that the pairing $\rho$ in Proposition \ref{dbiprod}
is invertible. Then according to Proposition \ref{twist-bialg} the
multiplication $\m_{C\otimes H\otimes B}^{\hat\rho}$ of the twisted
bialgebra $(C\otimes H\otimes B)_{\hat\rho}$ is given by
\begin{equation}\mathlabel{textmult-twist}
\m_{C\otimes H\otimes B}^{\hat\rho} =
\vstretch = 70
\begin{array}{c}
\divide\unitlens by 2
\begin{tangle}
\object{\sstyle C}\step[1.5]\object{\sstyle H}\step[3]\object{\sstyle B}%
 \step[4]\object{\sstyle C}\step[3]\object{\sstyle H}\step[1.5]%
 \object{\sstyle B}\\
\id\step\hcd\step[1.5]\cd\Step\cd\step[1.5]\hcd\step\id\\
\id\step\id\step\id\step\hdd\Step\x\Step\hd\step\id\step\id\step\id\\
\id\step\id\step\id\step\hrd\step\hddcd\Step\hdcd\step\hld\step%
 \id\step\id\step\id\\
\id\step\id\step\id\hstep\hdd\hstep\hx\step\id\Step\id\step\hx\hstep\hd%
\hstep\id\step\id\step\id\\
\id\step\id\step\id\hstep\Put(5,20)[cc]{\scriptstyle +}\hev\hstep%
 \hcd\hstep\hd\step\hdd\hstep\hcd\hstep\Put(5,20)[cc]{\scriptstyle -}\hev%
 \hstep\id\step\id\step\id\\
\id\step\id\step\id\Step\id\step\hx\step\hx\step\id\Step\id\step\id%
 \step\id\\
\id\step\id\step\d\step\lu\step\hcu\step\ru\step\dd\step\id\step\id\\
\id\step\id\Step\x\step[1.5]\id\step[1.5]\x\Step\id\step\id\\
\id\step\id\step\dd\Step\hd\step\id\step[1.5]\id\Step\d\step\id\step\id\\
\id\step\lu\step[3.5]\hcu\step[1.5]\id\step[3]\ru\step\id\\
\cu\step[4]\cu\step[3]\cu\\
\hstep{\sstyle C}\step[5.5]\object{\sstyle H}\step[5]\object{\sstyle B}%
\end{tangle}
\multiply\unitlens by 2
\\[-20pt]
\end{array}
\end{equation}
where the pairing $\rho$ is presented by
$\begin{array}{c}\\[8pt]
 \begin{tangle}\Put(5,20)[cc]{\scriptstyle +}\hev\end{tangle}
\end{array}$
and its convolution inverse $\rho^{\text{--}}$ by
$\begin{array}{c}\\[8pt]
\begin{tangle}\Put(5,20)[cc]{\scriptstyle -}\hev
\end{tangle}
\end{array}$.
\end{corollary}

\begin{proof}
Because of Proposition \ref{dbiprod} we can apply Definition
\ref{bialg-twist} and Proposition \ref{twist-bialg} to $\m_Z$ in
Proposition \ref{free-biprod}. A straightforward calculation then
yields the result.
\end{proof}
\abs
We will discuss in the following an example in which unavoidably
cross product bialgebras in certain braided categories emerge.

\begin{example}
{\normalfont Corollary \ref{double-biprod} is a (braided) generalization of
Majid's double biproduct construction \cite{Ma1:95}.
It has been shown in \cite{Ma1:95,Som1:96} that the quantum
enveloping algebra $\mathbf{U}$ in terms of Lusztig's construction
\cite{Lus1:93} is a double biproduct bialgebra. Explicitely,
let $(I,\cdot)$ be a Cartan datum, and 
$(X,Y,\langle .\,,.\rangle)$ be a root datum of type
$(I,\cdot)$. Given the commutative ring $k=\QQ(q)$, let $\mathbf{f}$ be
the $k$-algebra generated by $I$, factorized by the annihilator
radical of the unique
pairing $(.,.):k\langle I\rangle\times k\langle I\rangle\to k$
given in Proposition 1.2.3 in \cite{Lus1:93}. Furthermore let
$\mathbf{U}_0$ be the group algebra of $Y$ over $k$. Then \cite{Lus1:93},
$\mathbf{U}\cong \mathbf{f}\otimes\mathbf{U}_0\otimes\mathbf{f}$
are isomorphic Hopf algebras. The algebra $\mathbf{f}$ is both left
and right $\mathbf{U}_0$-crossed module bialgebra. The bialgebra
structure of $\mathbf{f}$ is induced by the algebra structure of 
$\mathbf{f}$ and by primitivity of all elements of
$k\langle I\rangle$. From \cite{Lus1:93} we know that
$\mathbf{f}=\bigoplus_{\nu\in\NN[I]} \mathbf{f}_\nu$ is an
$\NN[I]$-graded algebra. The root datum provides embeddings
$x:I\hookrightarrow X$, $i\mapsto x_i$ and $y:I\hookrightarrow Y$,
$i\mapsto y_i$ which in turn canonically induce homomorphisms of
abelian groups $x:\NN[I]\to X$, $\nu\mapsto x_\nu$ and $y:\NN[I]\to Y$,
$\nu\mapsto y_\nu$. Then the left $\mathbf{U}_0$-crossed module
structure of $\mathbf{f}$ is given by $y\triangleright \lambda :=
q^{-\langle y, x_{\vert\lambda\vert}\rangle}\,\lambda$ and 
$\nu_l(\lambda) := - \sum_{i\in I}(\frac{i\cdot i}2 \cdot\vert\lambda\vert_i)
\,y_i\otimes \lambda$ where $y\in Y$ and $\lambda\in f$ is
homogeneous of degree $\nu=\vert\lambda\vert=
\sum_{i\in I} \vert\lambda\vert_i\,i\in\NN[I]$, 
i.~e.~$\lambda\in\mathbf{f}_{\vert\lambda\vert}$. The right 
$\mathbf{U}_0$-crossed module structure of $\mathbf{f}$ is
$\lambda\triangleleft y := 
q^{-\langle y, x_{\vert\lambda\vert}\rangle}\,\lambda$ and 
$\nu_r(\lambda) := - \sum_{i\in I}(\frac{i\cdot i}2 \cdot\vert\lambda\vert_i)
\,\lambda\otimes y_i$. Then $\mathbf{U}\equiv \mathbf{f}\otimes
\mathbf{U}_0\otimes\mathbf{f}$ is a double biproduct bialgebra
according to Corollary \ref{double-biprod} with $H=\mathbf{U}_0$, 
$C=\mathbf{f}$ and $B=\mathbf{f}$, and the pairing $\rho=(.,.)$.

Although in the present example the base category $\C\,$$=$$\,k\text{-mod}$
is symmetric, the categories of $H$-crossed modules and $H$-Hopf bimodules
are braided and the cross product bialgebra construction is within
these categories. This emphasizes once again that the double biproduct
(even in ordinary symmetric categories) is a cross product bialgebra 
construction in a braided category.}
\end{example}

\subsection{Quantum Double Construction}

In Corollary \ref{hp-cases} we discussed braided versions of
matched pairs leading to braided double cross products. From
\cite{Ma2:94} we know that two dually paired bialgebras in a symmetric
category yield a matched pair from which a generalization of the
Drinfel'd double can be reconstructed. Such a procedure had been
discussed for braided categories in \cite{Ma1:95}. It was announced
that a similar construction fails there since
the braiding twists up and can not be disentangled. The subsequent
proposition confirms this observation in a certain sense.

\begin{proposition}\Label{braid-match}
Suppose that $A$ and $H$ are bialgebras in $\C$ which are paired by the
pairing $\langle .,.\rangle:H\otimes A\to \E_C$ subject to the defining
identities
\begin{equation}\mathlabel{match-pairing}
\begin{split}
\langle .,.\rangle\circ (\m_H\otimes\id) &=\langle .,.\rangle\circ
 (\id\otimes \langle .,.\rangle \otimes\id)\circ
 (\id\otimes\id\otimes\Psi\circ\Delta_A)\\
\langle .,.\rangle\circ (\id\otimes\m_A) &=
 \langle .,.\rangle\circ(\id\otimes\langle .,.\rangle\otimes\id)\circ
 (\Delta_H\otimes\id\otimes\id)\\
\langle .,.\rangle\circ(\eta_H\otimes\id) &= \varepsilon_A\\
\langle .,.\rangle\circ(\id\otimes\eta_A) &= \varepsilon_H\,.
\end{split}
\end{equation}
Then the following statements hold.
\begin{enumerate}
\item\Label{item1}
If $A$ and $H$ are Hopf algebras and $S_A$ is an isomorphism in $\C$ then
$\langle .,.\rangle$ is convolution invertible. Explicitely
$\langle .,.\rangle^{\text{--}}=\langle .,.\rangle\circ(S_H\otimes\id)=
\langle .,.\rangle\circ(\id\otimes S_A^{-1})$.
\item\Label{item2}
If $\langle .,.\rangle$ is convolution invertible
we define
\begin{equation}\mathlabel{match-pairing2}
\begin{split}
\lhd = &(\langle .,.\rangle^{\text{--}}\otimes\id\otimes\langle .,.\rangle)
             \circ(\id \otimes\Psi_{H\otimes H,A}\otimes\id)
             \circ(\Delta_H^{(2)}\otimes\Delta_A)\\
\rhd = &(\langle .,.\rangle^{\text{--}}\otimes\id\otimes\langle .,.\rangle)
             \circ(\id\otimes\Psi_{H,A\otimes A}\otimes\id)
             \circ(\Delta_H\otimes\Delta_A^{(2)})
\end{split}
\end{equation}
where $\Delta^{(2)}=(\Delta\otimes\id)\circ\Delta$. Then $(A,\rhd)$
is a left $H$-module and $(H,\lhd)$ is a right $A$-module.
The tuple $(A, H, \lhd,\rhd)$ is a matched pair as in Corollary
\ref{hp-cases} if and only if
$\Psi_{H,A}\circ\Psi_{A,H}=\id$.
\end{enumerate}
\end{proposition}

\begin{proof}
Statement \ref{item1} is proved analogously as in the standard
symmetric case. Without problems one verifies that $\lhd$ and $\rhd$
in statement \ref{item2} define algebra actions. Now suppose that
$\Psi_{H,A}\circ\Psi_{A,H}=\id_{A\otimes H}$. It is not difficult to show
that $(A, H, \lhd,\rhd)$ is a matched pair because nearly everything
works like in the classical symmetric case \cite{Ma1:90,Ma2:94}.
Conversely if $(A, H,\lhd,\rhd)$ defined by \eqref{match-pairing2} is a
matched pair then the following identity has to be fulfilled because of
the last equation in \eqref{braid-match-pair}.
\begin{equation}\mathlabel{braid-match1}
\vstretch = 70
\begin{array}{c}
\divide\unitlens by 2
\begin{tangle}
\step[2.5]\object{\sstyle H}\step[3]\object{\sstyle A}%
 \Step\hstep\\
\Step\hcd\Step\hcd\Step\\
\step\hddcd\step\x\step\hdcd\step\\
\dd\step\hx\Step\hx\step\d\\
\Pairing{\sstyle\langle,\rangle^{\text{--}}}\step\id\Step%
 \id\step\Pairing{\sstyle\langle,\rangle}\\
\step[3]\object{\sstyle H}\Step\object{\sstyle A}
\end{tangle}
\multiply\unitlens by 2
\end{array}
\quad = \quad
\begin{array}{c}
\divide\unitlens by 2
\begin{tangle}
\step\Step\object{\sstyle H}\Step\object{\sstyle A}%
 \Step\step\\
\step\Step\x\Step\step\\
\step\Step\x\Step\step\\
\Step\hddcd\Step\hdcd\Step\\
\step\hddcd\step\x\step\hdcd\step\\
\dd\step\hx\Step\hx\step\d\\
\Pairing{\sstyle\langle,\rangle^{\text{--}}}\step\id\Step
 \id\step\Pairing{\sstyle\langle,\rangle}\\
\step[3]\object{\sstyle H}\Step\object{\sstyle A}
\end{tangle}
\multiply\unitlens by 2
\end{array}
\end{equation}
From \eqref{braid-match1} we obtain $\Psi_{A,H}\circ\Psi_{H,A}=
\id_{H\otimes A}$ by multiplying $\langle .,.\rangle$ to the left
and $\langle .,.\rangle^{\text{--}}$ to the right of \eqref{braid-match1}
using the product given by \eqref{conv-prod}.
\end{proof}

\section{Conclusions and Outlook}\Label{conc-out}

We defined Hopf data and cross product bialgebras very generally.
Cross product bialgebras are Hopf data. A special class of Hopf
data are recursive Hopf data probably with finite order.
Recursive Hopf data are cross product bialgebras.
We further restricted to trivalent Hopf data and trivalent 
cross product bialgebra. We showed the equivalence of both notions and
provided a description of trivalent cross product bialgebras
either through (co-)modular properties or by universal systems
of certain projections and injections respectively.
Therefore the classification of trivalent cross product bialgebras
in terms of trivalent Hopf data has been achieved. The known cross
products with bialgebra structure fit into this new classification scheme.
In addition new types of trivalent cross product bialgebras have been found
which generalize all other types. However explicit examples have not
been found yet for these general types of trivalent cross product bialgebras.

We have been working throughout in a braided monoidal setting which allowed us
to apply the machinery of Hopf data and cross
product bialgebras to braided categories.
In particular we showed that the double biproduct bialgebras
come from a certain tensor product bialgebra in the braided category of
Hopf bimodules over a given Hopf algebra. A more general study of recursive 
Hopf data in Hopf bimodule categories will be published elsewhere.

The structure of (recursive) Hopf data shows to be symmetric under duality and
reflection at a vertical axis -- if one considers the defining identities
as graphics in three dimensional space. These symmetries will be
somehow destroyed when considering trivalent Hopf data and trivalent
cross product bialgebras and one might ask if such a breaking of
symmetry is a generic feature of the theory of cross product bialgebras.
Therefore it remains an open problem if the present setting is the most
general one to describe cross product bialgebras by (co-)modular
properties or by universal systems of projection and injection morphisms
equivalently. One could think of certain types of recursive Hopf data
(with finite order) or some other specializations of Hopf data which 
preserve the above mentioned symmetries, to be good candidates for a more 
general framework.
A possible generalization has been presented in Remark \ref{non-triv-coact}
although the symmetries will be destroyed in this case, too.

In the present article we studied cross product bialgebras
without cocycles and cycles (or dual cocycles). In a forthcoming
paper \cite{BD:98} we will apply similar techniques,
and results from \cite{BCM:86,Brz1:96,DT:86,Mon1:92} to describe
certain types of cross product bialgebras with co-cycles in a
co-cyclic (co-)modular way. 
Analogous statements as in Proposition \ref{crossprod-proj}
and Theorem \ref{rhp-cpb} will be derived
for rather general cross product bialgebras with co-cycles. But
the co-cyclic (co-)modular scheme of classification turns out to be
much more subtle than in the present case.
There might be different ways of restricting the general set-up
of co-cycle cross product bialgebra to achieve different sorts of
classification schemes.

\abs
\abs
\parbox{12cm}{\scriptsize{\textsc{Bogolyubov Institute for Theoretical
Physics, Metrologichna Str.~14-b, Kiev 252143, Ukraine.}
\textbf{E-mail:} {\tt mmtpitp@gluk.apc.org}
\\[0.5cm]
\textsc{DAMTP, University of Cambridge, Silver Street, Cambridge CB3 9EW,
UK.} \textbf{E-mail:} {\tt b.drabant@damtp.cam.ac.uk}}}

\end{document}